\documentclass[12pt]{article}
\usepackage{amsmath, setspace}
\usepackage[margin=1in]{geometry}
\usepackage[linesnumbered,ruled]{algorithm2e}
\usepackage{graphicx}
\usepackage[english]{babel}
\usepackage{amsthm}
\usepackage{amssymb}
\usepackage{array}
\usepackage{subfigure}
\usepackage[latin1]{inputenc}
\usepackage{enumerate}
\usepackage{xspace}
\usepackage{amsfonts}
\usepackage{amsmath}
\usepackage{hyperref}
\usepackage{caption}
\usepackage{booktabs}
\usepackage{color}
\usepackage{lmodern}
\usepackage{tikz}
\usepackage{comment}
\usepackage{enumitem}
\usepackage{pgfplots}
\usepackage{lipsum}
\usepackage{authblk}

\usepackage{color-edits}
\addauthor{zc}{cyan}
\addauthor{mc}{purple}
\addauthor{kt}{brown}

\usetikzlibrary{shapes}
\usetikzlibrary{snakes}
\usetikzlibrary{backgrounds}
\usetikzlibrary{calc}

\newtheorem{theorem}{Theorem}
\newtheorem{lemma}{Lemma}

\newtheorem{proposition}{Proposition}
\newtheorem{definition}{Definition}
\newtheorem{assumption}{Assumption}

\spacing{1.5}

\newcommand{\eg}{{\it e.g.}}
\newcommand{\ie}{{\it i.e.}}

\newcommand{\prom}{\textnormal{ProM}$^3$\xspace}

\DeclareMathOperator*{\argmin}{argmin}
\DeclareMathOperator*{\argmax}{argmax}

\DeclareMathOperator{\Proj}{Proj}

\newcommand{\bmt}[1]{\tilde{\bm{#1}}}

\usepackage{bm}
\usepackage{cleveref}
\usepackage{threeparttable}
\usepackage{ulem}

\usepackage{array}

\usepackage{natbib}
\bibpunct[, ]{(}{)}{,}{a}{}{,}%

\title{A Max-Min-Max Algorithm for \\Large-Scale Robust Optimization}

\author[1]{Kai Tu}
\author[2]{Zhi Chen}
\author[3]{Man-Chung Yue}

\affil[1]{\small \textit{Shenzhen University,} \texttt{kaitu\_02@163.com}}
\affil[2]{\small \textit{The Chinese University of Hong Kong,} \texttt{zhi.chen@cuhk.edu.hk}}
\affil[3]{\small \textit{The University of Hong Kong,} \texttt{mcyue@hku.hk}}
\date{\today}

\begin{document}
\maketitle

% \vspace{-28mm}
\begin{abstract}
Robust optimization (RO) is a powerful paradigm for decision making under uncertainty. Existing algorithms for solving RO, including the reformulation approach and the cutting-plane method, do not scale well, hindering the application of RO to large-scale decision problems. In this paper, we devise a first-order algorithm for solving RO based on a novel max-min-max perspective. Our algorithm operates directly on the model functions and sets through the subgradient and projection oracles, which enables the exploitation of problem structures and is especially suitable for large-scale RO. Theoretically, we prove that the oracle complexity of our algorithm for attaining an $\varepsilon$-approximate optimal solution is $\mathcal{O}(\varepsilon^{-3})$ or $\mathcal{O}(\varepsilon^{-2})$, depending on the smoothness of the model functions. The algorithm and its theoretical results are then extended to RO with projection-unfriendly uncertainty sets. We also show via extensive numerical experiments that the proposed algorithm outperforms the reformulation approach, the cutting-plane method and two other recent first-order algorithms.

\end{abstract}

\section{Introduction}\label{sec:introduction}

Optimization models often require the input of some instance-specific parameters, which are unfortunately uncertain in most applications. The uncertainty could come from errors in estimating or measuring the parameters. Another reason for the uncertainty could be that the parameters are intrinsically random. For example, the beamforming problem in wireless communication (\citealt{wu2017sdr}) aims at finding the optimal transmission angle and power concerning some objective (\eg, least interference or largest throughput) subject to certain physical constraints. The parameters required to specify the beamforming optimization model include the transmitter and receiver antennas' geographical locations, the obstacles between them, and the spectrum of other network users.  Misspecification of these parameters may lead to poor signal quality or even a breakdown of the communication network. As one of the most powerful and popular paradigms for optimization under uncertainty, robust optimization (RO) has attracted intense research in recent years and found applications across a wide range of areas such as machine learning~(\citealt{singla2020survey}), operations management~(\citealt{bertsimas2023data}), health care~(\citealt{meng2015robust}) and finance~(\citealt{gregory2011robust}), to name a few.

To set the scene, consider the following nominal optimization problem:
\[
\begin{array}{rll}
\min & f_0(\bm{x}) \\
{\rm s.t.} & g_m(\bm{x}, \bm{z}_m) \leq 0 &~\forall m \in [M] \\
& \bm{x} \in \mathcal{X},
\end{array}
\]
where $f_0$ is the objective function, $g_1,\dots, g_M$ are the constraint functions, $\bm{x}\in\mathbb{R}^N$ is the decision vector, the set $\mathcal{X}\subseteq \mathbb{R}^N$ models further constraints on $\bm{x}$, and $\bm{z}_1 \dots, \bm{z}_M\in\mathbb{R}^{J_m}$ are the parameters specifying the optimization model. In the face of uncertainty, RO postulates that each parameter $\bm{z}_m $ resides in a subset $\mathcal{Z}_m \subseteq \mathbb{R}^{J_m} $---called the uncertainty set---that represents the modeler's belief about the possible range of the uncertain parameter $\bm{z}_m$. RO then takes a pessimistic point of view: whatever decision is chosen, the worst parameters over the uncertainty sets will be realized accordingly. More precisely, RO prescribes choosing the decision as an optimal solution to the problem
\begin{equation}\label{opt:RO1}\tag{\sc Robust}
\begin{array}{rll}
\min & f_0(\bm{x}) \\
{\rm s.t.} & \displaystyle \max_{\bm{z}_m \in \mathcal{Z}_m} g_m (\bm{x}, \bm{z}_m) \leq 0  \quad\forall m\in [M] \\
& \bm{x} \in \mathcal{X}.
\end{array}
\end{equation}
Despite the wide applicability, computational approaches for solving RO are to some extent limited and do not meet the needs of modern RO users who often require to solve high-dimensional non-linear and/or non-smooth RO problems. 
The difficulty stems primarily from the embedded optimization problems $\max_{\bm{z}_m\in\mathcal{Z}_m} g_m(\bm{x},\bm{z}_m)$, rendering standard optimization algorithms non-viable. 

Currently, there are two dominant approaches: the reformulation approach (see, \eg, \citealt{ben1998robust}) and the cutting-plane method (see, \eg, \citealt{mutapcic2009cutting}). Roughly speaking, the reformulation approach solves the \ref{opt:RO1} problem by converting it into a deterministic reformulation---an equivalent optimization problem without embedded optimization problems. Here, ``deterministic'' refers to the uncertainty-free nature of the reformulation. This is often achieved by either solving the embedded optimization problem analytically or replacing it with its dual problem. Since the deterministic reformulation is a standard optimization problem, it can be solved by many sophisticated off-the-shelf optimization solvers such as CPLEX, Gurobi and MOSEK. The reformulation approach works for a large and useful class of \ref{opt:RO1} problems wherein the constraint functions $g_1, \dots, g_M$ and the uncertainty sets $\mathcal{Z}_1,\dots, \mathcal{Z}_M$ possess certain special structures; see~\cite{ben1998robust, ben2009robust}.

The cutting-plane method is essentially Kelley's cutting-plane method (\citealt{kelley1960cutting}) specialized to RO. It is an iterative algorithm alternating between two steps: the \textit{optimization step} and the \textit{pessimization step}. Given finite subsets $\widehat{\mathcal{Z}}_m \subseteq \mathcal{Z}_m$, $m\in[M]$, the optimization step solves the following approximation of \ref{opt:RO1} to find a new $\bm{x}$:
\begin{equation*}\label{opt:RO1_discrete}
\begin{array}{rl}
\min & f_0(\bm{x}) \\
{\rm s.t.} & \displaystyle \max_{\bm{z}_m \in\widehat{\mathcal{Z}}_m } g_m(\bm{x},\bm{z}_m) \le 0 \quad\forall m\in [M] \\
& \bm{x}\in\mathcal{X}.
\end{array}
\end{equation*}
Such an approximation has only a finite number of constraints and thus can be solved readily. The pessimization step computes a maximizer $\hat{\bm{z}}_m$ for each embedded problem $\max_{\bm{z}_m\in\mathcal{Z}_m} g_m(\bm{x},\bm{z}_m)$ and if the maximum value is positive, expands the approximate uncertainty set $\widehat{\mathcal{Z}}_m$ by setting $\widehat{\mathcal{Z}}_m \leftarrow \widehat{\mathcal{Z}}_m\cup \{ \hat{\bm{z}}_m \}$. The algorithm alternates between these two steps until reaching a certain stopping criterion.

There are several weaknesses of existing computational approaches to RO, which significantly hinder its development and applications. 
First, off-the-shelf optimization solvers typically applied to solve the deterministic reformulations of RO are mainly based on interior-point methods (\citealt{nesterov1994interior}). As evidenced by many numerical studies~(\citealt{toh2010accelerated, liu2023unified}), interior-point methods have relatively worse scalability when compared to some other classes of optimization algorithms, such as first-order and second-order methods, and thus are unsuitable for large-scale problems. Moreover, solvers may ignore useful structures of the problem (\eg, sparsity or low-rank-ness of matrices), which, if suitably exploited, can substantially improve the computational speed. 

Second, as mentioned earlier, one way to eliminate the embedded maximization in~\ref{opt:RO1} is to dualize it into a minimization problem. The dualization process possibly introduces a large number of extra variables and/or constraints, resulting in a high-dimensional and/or highly constrained reformulation that is hard to solve even for well-developed solvers.

Third, for many cases of $g_m$ and $\mathcal{Z}_m$, the reformulation technique often lifts \ref{opt:RO1} to a more difficult and general class of optimization problems. As an example, suppose that the constraint function is
$ g_m (\bm{x},\bm{z}_m) = \bm{x}^\top\bm{A}^\top(\bm{z}_m)\bm{A}(\bm{z}_m)\bm{x} + \bm{b}^\top(\bm{z}_m)\bm{x} + c(\bm{z}_m) 
$,
with $\bm{A}(\bm{z}_m)\in \mathbb{R}^{L \times N}$, $\bm{b}(\bm{z}_m)\in\mathbb{R}^N$ and $c(\bm{z}_m)\in\mathbb{R}$ being affine functions in $\bm{z}_m$, and the uncertainty set $\mathcal{Z}_m$ is a unit ball in $\mathbb{R}^{J_m}$.
In this case, the function $g_m$ is quadratic in both $\bm{x}$ and $\bm{z}_m$, and the uncertainty set $\mathcal{Z}_m$ is defined by a quadratic inequality. Its reformulation, however, is a semidefinite programming problem~(\citealt{ben1998robust}). The quadratic nature of the robust constraint is disregarded.

Fourth, it is well-known that the cutting-plane method may perform poorly in both theory and practice (\citealt{mitchell2009cutting}) due to instability. In particular, its iteration complexity (\ie, the number of iterations required to achieve an $\varepsilon$-optimal solution) is $(1+ \mathcal{O}(\varepsilon^{-1}))^N$, exponential in the dimension~\cite[Section~5.2]{mutapcic2009cutting}. Therefore, its performance on RO problems could potentially be equally bad. Moreover, it has been shown in a comprehensive computational study~(\citealt{bertsimas2016reformulation}) that the cutting-plane method performs on par with or even worse than the reformulation approach.

Motivated by the above discussions, we aim to develop specialized iterative algorithms for efficiently solving large-scale \ref{opt:RO1} problems. The idea of specialized iterative algorithms for RO is not entirely new and has been recently pursued. Using tools from online convex optimization (\citealt{hazan2019introduction}), \cite{ben2015oracle} developed an iterative algorithm with an iteration complexity $\mathcal{O}(\varepsilon^{-2})$ for solving RO, each iteration of which requires an optimization step similar to that in the cutting-plane method.
Extending the online convex optimization idea, significant improvements have been obtained in the papers~\citealt{ho2018online, ho2019exploiting}, where the authors developed a first-order method (each iteration requires only first-order updates but neither the optimization nor pessimization step) with an iteration complexity of $\mathcal{O}(\varepsilon^{-1}\log \tfrac{1}{\varepsilon})$. One drawback of this algorithm is that it requires a binary search of the optimal value, which incurs extra computational overhead. In the very recent work~\cite{postek2021first}, another first-order method, named SGSP, has been developed based on perspective transformations~(see \citealt{boyd2004convex} for reference). Its iteration complexity is $\mathcal{O}(\varepsilon^{-2})$, with respect to more complicated oracles than the ones assumed in~\citealt{ho2018online, ho2019exploiting} and this paper tough. 
% Experiment results presented in~\cite{postek2021first} suggested that the numerical performance of both the cutting-plane method and the algorithm in~\cite{ho2018online, ho2019exploiting} are \zccomment{inferior} to SGSP.

The point of departure of our work is the following max-min-max problem:
\begin{equation}\label{opt:mmm1}\tag{{\sc Max-Min-Max}}
\max_{\bm{\lambda}\ge \bm{0}} \; \min_{\bm{x}\in\mathcal{X}} \; \max_{\bm{z} \in \mathcal{Z}} \; \mathcal{K}(\bm{\lambda}, \bm{x}, \bm{z})  ,
\end{equation}
where $\mathcal{K}(\bm{\lambda}, \bm{x}, \bm{z}) = f_0(\bm{x}) + \sum_{m \in [M]} \lambda_{m} g_m(\bm{x},\bm{z}_m) $ and $\mathcal{Z} = \mathcal{Z}_1\times\cdots\times\mathcal{Z}_M$.
This is essentially the Lagrangian dual of the \ref{opt:RO1} problem and therefore equivalent to it under mild assumptions (see \Cref{prop:mmm-equivalence} below).
The advantage of our framework is that neither complicated reformulation (as in the reformulation approach) nor perspective transformation (as in \citealt{postek2021first}) is needed. Furthermore, the algorithm will operate directly on the functions $f_0, g_1, \dots, g_M$ as well as the sets $\mathcal{X},\mathcal{Z}_1, \dots, \mathcal{Z}_M$ through their gradient and projection oracles, respectively. The useful structures and theoretical properties of the functions and sets are all preserved and can be easily exploited in algorithmic design.
However, designing and theoretically analyzing max-min-max algorithms are substantially more difficult than those for max-min problems. Our contributions are as follows.
\begin{itemize}
\item We devise a first-order method, called \prom, for solving the \ref{opt:mmm1} problem (hence, the \ref{opt:RO1} problem) that utilizes the structure of $\mathcal{K}(\bm{\lambda}, \bm{x}, \bm{z})$ and operates directly on the constituent functions and sets. 
Our design views~\ref{opt:mmm1} as a max-min problem 
\begin{equation}\label{opt:outer-saddle}
\max_{\bm{\lambda}\ge \bm{0}} \; \min_{\bm{x}\in\mathcal{X}} \; \mathcal{L}(\bm{\lambda}, \bm{x}),
\end{equation}
where $\mathcal{L}(\bm{\lambda}, \bm{x}) = \max\limits_{\bm{z} \in \mathcal{Z}} \mathcal{K}(\bm{\lambda}, \bm{x}, \bm{z})$ is the Lagrangian function, and adopts the framework of alternating proximal algorithm for max-min problems. We call problem~\eqref{opt:outer-saddle} the \textit{outer saddle-point problem}.
Nevertheless, since the objective function~$\mathcal{L}(\bm{\lambda}, \bm{x})$ in the outer saddle-point problem is itself a maximum value, efficiently updating $\bm{\lambda}$ and $\bm{x}$ as prescribed by the standard alternating proximal algorithm is non-trivial.
Our algorithm therefore requires extra ideas. In particular, the updating step for $\bm{x}$ turns out to be a strongly-convex-concave min-max problem with a non-linear coupling term, which we call the \textit{inner saddle-point problem}. A customized algorithm for the inner saddle-point problem is also developed. Combining the algorithms proposed for the outer and inner saddle-point problems leads to our algorithm \prom for solving RO.

\item We prove that under similar conditions as in \cite{ho2018online, ho2019exploiting, postek2021first}, the proposed algorithms for both the outer and inner saddle-point problems enjoy a sublinear convergence rate. 
Since the per-iteration cost of different algorithms may vary, the oracle complexity---the number of calls to the projection and subgradient oracles---for achieving an $\varepsilon$-optimal solution would be a more faithful measure of computational efficiency. 
Based on the convergence analysis of our algorithms for the outer and inner saddle-point problems, we prove that our algorithm \prom enjoys the oracle complexity of $\mathcal{O}(\varepsilon^{-3})$, and the oracle complexity can be strengthened to $\mathcal{O}(\varepsilon^{-2})$ if the functions $f$ and $g_1,\dots,g_M$ are all smooth.

\item Third, we extend our algorithm \prom and convergence analysis to RO problems where the uncertainty sets $\mathcal{Z}_1,\dots,\mathcal{Z}_M$ take certain intersection form and do not admit easy projection. Such a setting is particularly useful in distributionally robust optimization, where the uncertainty $\bm{z}_m$ is a probability vector lying in the intersection of the probability simplex and a ball defined by some norm, such as those based on the popular (type-$\infty$) Wasserstein distance (\citealt{Esfahani_Kuhn_2017, xie2019tractable, bertsimas2022two, gao2023distributionally}) or the support space is a finite set (\citealt{Ben_Hertog_Waegenaere_Melenberg_Rennen_2013, wiesemann2014distributionally}). 
\end{itemize}

We conclude the introduction with a few remarks. First, although the oracle complexity of SGSP in~\cite{postek2021first} is $\mathcal{O}(\varepsilon^{-2})$, as mentioned above, it relies on oracles that are generally more complicated than those assumed in this paper. Thus, the complexity $\mathcal{O}(\varepsilon^{-2})$ or $\mathcal{O}(\varepsilon^{-3})$ of \prom should not be directly compared to the complexity $\mathcal{O}(\varepsilon^{-2})$ of SGSP. 
Second, contrary to saddle-point problems that have received intense research recently, the literature on max-min-max problems, as pointed out by~\citealt{polak2003algorithms}, is much scarcer. Therefore, our work could be of independent interest to researchers working on max-min-max problems. 
Third, to the best of our knowledge, strongly-convex-concave min-max problems with a non-linear coupling term have not been thoroughly investigated yet. Our proposed algorithm for the inner saddle-point problem and its convergence analysis partially fill this gap in the fast-growing literature on saddle-point problems. 

\section{Proximal Max-Min-Max Algorithm (ProM$^3$)}\label{sec:prom}

We develop a first-order algorithm for solving the~\ref{opt:mmm1} problem (thus, the \ref{opt:RO1} problem), called the proximal max-min-max algorithm (\prom). Our line of attack to~\ref{opt:mmm1} is to view it as two layers of saddle-point problems: the outer saddle-point problem is the max-min problem~\eqref{opt:outer-saddle}, whereas the inner saddle-point problem is a step towards solving the outer problem in our algorithmic framework. For simplicity, we call the algorithms for solving the outer and inner saddle-point problems the outer and inner algorithms, respectively. The proposed \prom for solving~\ref{opt:mmm1} is obtained by combining the outer and inner algorithms. 

\subsection{Outer Algorithm}\label{sec:outer_alg}

To describe the outer algorithm, we re-state the outer saddle-point problem~\eqref{opt:outer-saddle} as
\begin{equation}\label{opt:outer-saddle-1}\tag{\sc Outer}
\max_{\bm{\lambda}\ge \bm{0}} \; \min_{\bm{x}\in\mathcal{X}} \;  
\left\{ 
f_0(\bm{x}) + \bm{\lambda}^\top \bm{f}(\bm{x})
\right\},
\end{equation}
where $\bm{f}(\bm{x}) = (f_1(\bm{x}), \dots, f_M(\bm{x}))$ and $f_m(\bm{x}) = \max_{\bm{z}_m\in\mathcal{Z}_m} g_m(\bm{x},\bm{z}_m)$ for each $m\in [M]$.
The general algorithmic framework we adopt for the outer saddle-point problem, \ref{opt:outer-saddle-1}, is an alternating proximal update scheme: 
\begin{equation}\label{label: Original-PD}
\left\{
\begin{array}{r@{\;}l@{\;}l}
\bm{\lambda}^{k+1} & = & \displaystyle \argmax_{\bm{\lambda} \geq \bm{0}} \; \bm{\lambda}^\top \bm{f}(\bm{x}^k) - \frac{1}{2\beta}\|\bm{\lambda} - \bm{\lambda}^k\|_2^2 =  [\bm{\lambda}^k + \beta \bm{f}(\bm{x}^k)]_+,\\
\bm{x}^{k+1} & = & \displaystyle \argmin_{\bm{x} \in \mathcal{X}} \; f_0(\bm{x})+ (\bm{\lambda}^{k+1})^\top  \bm{f}( \bm{x}) + \frac{1}{2\alpha} \|\bm{x} - \bm{x}^k \|_2^2, \\
\end{array}
\right.
\end{equation}
where $\alpha,\beta >0$ are step sizes and $[a]_+ = \max\{a, 0\}$ for any $a\in\mathbb{R}$. It is well-known that scheme~\eqref{label: Original-PD} is generally divergent. A correction for the scheme is proposed in~\cite{chambolle2011first}, which has a modified $\bm{\lambda}$-update but keeps the $\bm{x}$-update unchanged. When specialized to our \ref{opt:outer-saddle-1} problem, the modified $\bm{\lambda}$-update reads 
\[
\bm{\lambda}^{k+1} = \displaystyle \argmax_{\bm{\lambda} \geq \bm{0}} \; \bm{\lambda}^\top (2 \bm{f}(\bm{x}^k) - \bm{f}(\bm{x}^{k-1})) - \frac{1}{2\beta}\|\bm{\lambda} - \bm{\lambda}^k \|_2^2 = [\bm{\lambda}^k + \beta (2 \bm{f}(\bm{x}^k) - \bm{f}(\bm{x}^{k-1}))]_+.
\]
Unfortunately, even the corrected scheme would not work in our situation, at least not efficiently. 
The reason is that each component $f_m (\bm{x}^k)$ of the vector $\bm{f}(\bm{x}^k)$ is a partial maximum of $g_m (\bm{x}^k, \bm{z}_m)$ with respect to its second argument $\bm{z}_m$. This prohibits efficient evaluation of $f_m$ or its gradient. We circumvent this issue by approximating $\bm{f}$ in both $\bm{\lambda}$- and $\bm{x}$-updates. 
Specifically, to approximate the $\bm{\lambda}$-update, we replace $\bm{f}(\bm{x}^k)$ by
\[ 
\bm{g}(\bm{x}^k, \bm{z}^k) = (g_1( \bm{x}^k, \bm{z}_1^k ),\dots, g_M( \bm{x}^k, \bm{z}_M^k )),
\]
where each $ \bm{z}_m^k $ is an approximate maximizer of $g_m( \bm{x}^k, \cdot )$ over $\mathcal{Z}_m$  satisfying the condition
\[  
f_m (\bm{x}^k) = \max_{\bm{z}_m \in \mathcal{Z}_m} \; g_{m}(\bm{x}^{k},~ \bm{z}_{m}) \le g_{m}(\bm{x}^{k},~ \bm{z}_{m}^k) +\theta
\]
for some prescribed $\theta > 0$. 

For the $\bm{x}$-update, by noting that it is equivalent to the min-max problem
\[
\min_{\bm{x} \in \mathcal{X}} \; \max_{\bm{z} \in \mathcal{Z}} \; 
\left\{
f_0(\bm{x}) + \sum_{m\in [M] }\lambda_{m}^{k+1} g_m(\bm{x}, \bm{z}_m) + \frac{1}{2\alpha}\|\bm{x} - \bm{x}^k \|_2^2
\right\},
\]
we invoke a min-max algorithm to solve it to a custom-made stopping condition below.

\begin{definition}[Strong Approximate Saddle Point]\label{def: New-approximation-SP}
Consider a function $\mathcal{F}:\mathcal{U}\times \mathcal{V} \to \mathbb{R}$ such that $\mathcal{F}(\cdot, \bm{v})$ is $\sigma$-strongly convex on $\mathcal{U}$ for any $\bm{v}\in \mathcal{V}$, where $\sigma > 0$ is some constant independent of $\bm{v}$, and $\mathcal{F}(\bm{u}, \cdot)$ is concave on $\mathcal{V}$ for any $\bm{u}\in \mathcal{U}$. For any $\nu >0$, a pair $(\bmt{u}, \bmt{v})\in \mathcal{U}\times \mathcal{V}$ is said to be a strong $\nu$-approximate saddle point of $\mathcal{F}$ if
\[
\mathcal{F}(\bmt{u}, \bm{v}) \le \mathcal{F}( \bm{u}, \bmt{v} ) - \frac{\sigma}{2} \| \bm{u} - \bmt{u} \|_2^2 +\nu ~~~\forall (\bm{u},\bm{v})\in \mathcal{U} \times \mathcal{V}.
\]
\end{definition}

\Cref{def: New-approximation-SP} is stronger than the standard notion of approximate saddle point concept, $\mathcal{F}(\bmt{u}, \bm{v}) - \mathcal{F}( \bm{u}, \bmt{v} ) \le \nu$. The following proposition guarantees the existence of strong approximate saddle points under mild assumptions.

\begin{proposition}\label{prop:str_approx_sad_pt_exist}
Let $\mathcal{U}$ and $\mathcal{V}$ be  non-empty compact convex sets and $\mathcal{F}:\mathcal{U}\times \mathcal{V}\to \mathbb{R}$ be a function such that $\mathcal{F}(\cdot, \bm{v})$ is $\sigma$-strongly convex on $\mathcal{U}$ for any $\bm{v}\in \mathcal{V}$, where $\sigma > 0$ is some constant independent of $\bm{v}$, and $\mathcal{F}(\bm{u}, \cdot)$ is concave on $\mathcal{V}$ for any $\bm{u}\in \mathcal{U}$. Then, for any $\nu > 0$, $\mathcal{F}$ has a strong $\nu$-approximate saddle point.
\end{proposition}

The outer algorithm is formally presented in Algorithm~\ref{alg:outer}, and its convergence analysis will be presented in Section~\ref{subsec:out_rate}. It should be pointed out that although the alternating proximal algorithm for saddle-point problems has been studied, the convergence rate for its outer algorithm does not readily follow from existing results but requires certain new ideas; see the discussion after \Cref{thm: iteration-complexity-Algo1} for details.

\begin{algorithm}[t]
\caption{Outer Algorithm. \label{alg:outer}}
\SetKwInOut{Input}{Input}
\SetKwInOut{Output}{Output}
\Input{$K\geq 1$, $\theta>0$, $\nu> 0$, $\alpha>0$, $\beta>0$,  $\bm{\lambda}^{0} = \bm{0}$ and $ \bm{x}^{0}\in \mathcal{X}$.}
\For{$k= 0,1, \dots, K-1$}
{
($\bm{\lambda}$-update) For $m\in [M]$, find $\bm{z}_{m}^{k} \in \mathcal{Z}_m$ satisfying $f_{m}(\bm{x}^{k}) - g_{m}(\bm{x}^{k},~ \bm{z}_{m}^{k}) \leq \theta$. Set  
\[
\bm{\lambda}^{k+1} = [\bm{\lambda}^k + \beta (2\, \bm{g}(\bm{x}^k, \bm{z}^k) - \bm{g}(\bm{x}^{k-1}, \bm{z}^{k-1}))]_+ ,
\]
where $\bm{x}^{-1} = \bm{x}^0$ and $\bm{z}^{-1} = \bm{z}^0$ for the $0$-th iteration. \vspace{3mm}

($\bm{x}$-update) Compute a strong $\nu$-approximate saddle point $(\bm{x}^{k+1}, \bmt{z}^{k+1})\in \mathcal{X}\times \mathcal{Z}$ of the problem
\begin{equation}\label{opt:inner-saddle}\tag{\sc Inner}
\min_{\bm{x} \in \mathcal{X}} \; \max_{\bm{z} \in \mathcal{Z}} \; 
\left\{
f_0(\bm{x}) + \sum_{m\in [M] }\lambda_{m}^{k+1} g_m(\bm{x}, \bm{z}_m) + \frac{1}{2\alpha} \|\bm{x}- \bm{x}^k\|_2^2
\right\}.
\end{equation}
}
\Output{$\bar{\bm{x}}^K =\frac{1}{K} \sum_{k\in [K] }\bm{x}^{k}$.}
\end{algorithm}

Here, we offer an observation that can improve the practical performance of the outer algorithm. In the $\bm{\lambda}$-update, we do not necessarily need to compute the approximate maximizer $\bm{z}_m^k$ for every $m\in[M]$. Indeed, if $\lambda_m^k > 0$, then by definition, $\bmt{z}_m^k$ approximately maximizes the function $g_m(\bm{x}^k, \cdot )$ over $\mathcal{Z}_m$ and hence, with a careful choice of $\theta$ and $\nu$, is precisely the $\bm{z}_m^k$ that we need to compute at the beginning of the next iteration. 

\subsection{Inner Algorithm}\label{sec:inner_alg}

The $\bm{x}$-update step in the outer algorithm (see step~3 of Algorithm~\ref{alg:outer}) is the inner saddle-point problem, \ref{opt:inner-saddle}. 
A distinctive property of \ref{opt:inner-saddle} is a non-linear and non-smooth coupling term between the two variables $\bm{x}$ and $\bm{z}$. This should be contrasted with the relatively more common assumption that the coupling term is bilinear, \ie, $\bm{x}^\top \bm{Q} \bm{z}$ for some matrix $\bm{Q}$. Another feature of the \ref{opt:inner-saddle} problem is that its objective function is strongly convex in $\bm{x}$.
Saddle-point problems of this specific form have not been well studied. 

Due to the generality of the coupling term, algorithmic options are limited. We adopt the following variant of the subgradient ascent descent algorithm, which is known to enjoy a better convergence rate in the smooth case:
\[
\bm{z}_{t+1} = \Proj_{\mathcal{Z}} \left( \bm{z}_t - \delta ( 2\bm{\zeta}_t - \bm{\zeta}_{t-1} )\right),
\]
where $\bm{\zeta}_t \in \partial_{\bm{z}}(-f_0 - (\bm{\lambda}^{k+1})^\top \bm{g})(\bm{x}_t, \bm{z}_t)$ and $\Proj_{\mathcal{U}}(\cdot)$ denotes the projection onto a close convex set $\mathcal{U}$.
As the maximization over $\bm{z}$ in the \ref{opt:inner-saddle} problem is decomposable, the modified subgradient step can be executed by updating each $\bm{z}_m$ separately: that is, for each $m\in [M]$,
\begin{equation*}
\bm{z}_{t+1,m} = \Proj_{\mathcal{Z}_m} \left( \bm{z}_{t,m} - \delta \,\lambda_m^{k+1}\, (2\bm{\zeta}_{t,m} - \bm{\zeta}_{t-1,m} ) \right),
\end{equation*}
where $\bm{\zeta}_{t,m} \in \partial_{\bm{z}_m} (-g_m)( \bm{x}_t, \bm{z}_{t,m} )$. 

For the $\bm{x}$-update (the subgradient descent step), we slightly tweak the standard subgradient ascent descent framework to exploit the strong convexity.
Specifically, in the $\bm{x}$-update, we linearize only the non-smooth part but retain the strongly convex quadratic term: 
\[
\begin{array}{r@{\;}l@{\;}l}
\bm{x}_{t+1} & = & \displaystyle \argmin_{ \bm{x}\in\mathcal{X} }\; 
\left\{ 
\bm{\xi}_t^\top (\bm{x} - \bm{x}_t) + \frac{1}{2\alpha} \| \bm{x} - \bm{x}^{k-1} \|_2^2 + \frac{1}{2\gamma} \| \bm{x} - \bm{x}_t \|_2^2 
\right\} \\[3mm]
& = & \displaystyle \Proj_{\mathcal{X}} \left( \frac{\alpha \gamma}{\alpha + \gamma}\left(\frac{1}{\alpha}\bm{x}^{k-1} + \frac{1}{\gamma} \bm{x}_t - \bm{\xi}_t\right) \right),
\end{array}
\]
where $\bm{\xi}_t \in \partial_{\bm{x}} (f_0(\bm{x}_t) + (\bm{\lambda}^{k+1})^\top \bm{g}(\bm{x}_t, \bm{z}_{t+1}) )$.
This tweak allows for a larger step size $\gamma$ and improves practical performance. 
Note that by the subdifferential sum rule, the desired subgradient $\bm{\xi}_t$ can be obtained by $\bm{\xi}_t = \bm{\xi}_{t,0} + \lambda_1^{k+1}\, \bm{\xi}_{t,1}+\cdots+ \lambda_M^{k+1}\, \bm{\xi}_{t,M}$,
where $\bm{\xi}_{t,0} \in \partial f_0 (\bm{x}_t)$ and $\bm{\xi}_{t,m} \in \partial_{\bm{x}} g_m (\bm{x}_t, \bm{z}_{t+1,m})$ for $m\in[M]$. 
The inner algorithm is formally presented in Algorithm~\ref{alg:inner}, and its convergence analysis will be presented in Section~\ref{subsec:in_rate}.
\begin{algorithm}[t]
\caption{Inner Algorithm. \label{alg:inner}}
\SetKwInOut{Input}{Input}
\SetKwInOut{Output}{Output}
\Input{$T\ge 1$,  $\delta >0$,  $\gamma>0$, $\bm{x}_{0} \in \mathcal{X}$, $\bm{z}_0\in \mathcal{Z}$, $\alpha>0$, $\bm{x}^k\in \mathcal{X}$ and $\bm{\lambda}^{k+1} \ge \bm{0}$.}
\For{$t = 0, \dots, T-1$}
{
($\bm{z}$-update) For $m\in [M]$, compute $\bm{\zeta}_{t,m} \in \partial_{\bm{z}_m} (-g_m)( \bm{x}_t, \bm{z}_{t,m} )$. Set
\[ 
\bm{z}_{t+1,m} = \Proj_{\mathcal{Z}_m} \left( \bm{z}_{t,m} - \delta \,\lambda_m^{k+1}\, (2\bm{\zeta}_{t,m} - \bm{\zeta}_{t-1,m} ) \right), 
\]
where $\bm{\zeta}_{-1,m} = \bm{\zeta}_{0,m}$ for the $0$-th iteration. \vspace{3mm}

($\bm{x}$-update) For $m\in[M]$, compute $\bm{\xi}_{t,0} \in \partial f_0 (\bm{x}_t)$ and $\bm{\xi}_{t,m} \in \partial_{\bm{x}} g_m (\bm{x}_t, \bm{z}_{t+1,m})$. Set $\bm{\xi}_t = \bm{\xi}_{t,0} + \lambda_1^{k+1}\, \bm{\xi}_{t,1}+\cdots+ \lambda_M^{k+1}\, \bm{\xi}_{t,M}$ and
\[
\bm{x}_{t+1}  =  \Proj_{\mathcal{X}} \left( \frac{\alpha \gamma}{\alpha + \gamma}\left(\frac{1}{\alpha}\bm{x}^k + \frac{1}{\gamma} \bm{x}_t - \bm{\xi}_t\right) \right).
\]
}
\Output{$\bar{\bm{x}}_T =\frac{1}{T} \sum_{t\in [T]}\bm{x}_{t} $ and $\bar{\bm{z} }_T = \frac{1}{T}\sum_{t\in [T]}\bm{z}_{t}$.}
\end{algorithm}

Although the inner algorithm and its theoretical results are developed and presented in relation to~\ref{opt:inner-saddle}, they are actually applicable to general saddle-point problems of the form
\[  
\min_{\bm{u}\in \mathcal{U}} \; \max_{\bm{v}\in \mathcal{V}} \; \mathcal{F}(\bm{u}, \bm{v}),
\]
where the objective function $\mathcal{F}(\bm{u}, \bm{v})$ is strongly convex in $\bm{u}$, concave in $\bm{v}$ and has a general non-linear and non-smooth coupling term, and where the feasible regions $\mathcal{U}$ and $\mathcal{V}$ are non-empty closed convex sets. 
In fact, our proofs for the convergence results are presented in this general setting; see Appendix~\ref{sec:inner-conv-analysis}. 
However, for the convenience of RO theorists and practitioners, we customize the presentation of the inner algorithm to \ref{opt:inner-saddle} in the main text.

We also remark that our proposed algorithm \prom for RO relies only on the projection and subgradient oracles for the sets and functions, respectively, that define the \ref{opt:RO1} problem, and can be easily implemented.
Contrary to SGSP~(\citealt{postek2021first}), we do not need to pre-compute a Slater point or an upper bound of the dual optimal solution.
Numerical experiments in Section~\ref{sec:exp} show its promising performance on large-scale instances, in comparison with the reformulation approach, cutting-plane method, and the first-order methods developed in \cite{postek2021first} and \cite{ho2018online}. 

\section{Convergence Analysis}\label{sec:convergence_analysis}

This section determines the oracle complexity of the proposed algorithm \prom. To this end, we first theoretically analyze the convergence behavior of the outer and inner algorithms.

We collect the assumptions needed for our theoretical development. 
Assumption~\ref{Assump: Robust-Prob-1} below is standard in the RO literature (see, \eg, \citealt{ben2009robust}).

\begin{assumption}\label{Assump: Robust-Prob-1}
The following conditions hold.
\begin{enumerate}[font=\normalfont, label=(\textit{\roman*})]
\item\label{Assump: Robust-Prob-1-i} {\normalfont (Compactness and Convexity of Sets)} The sets $\mathcal{Z}_1\subseteq \mathbb{R}^{J_1},\dots, \mathcal{Z}_M\subseteq \mathbb{R}^{J_M}$ and $\mathcal{X}\subseteq \mathbb{R}^N$ are non-empty, compact and convex.
		
\item\label{Assump: Robust-Prob-1-ii} {\normalfont (Convexity of Functions)} 
The function $f_0: \mathbb{R}^N \to \mathbb{R}\cup\{+\infty\}$ is convex, with $\mathcal{X}\subseteq \mathrm{dom}(f_0)$. 
For any $m\in[M]$, the function $g_m:\mathbb{R}^N\times \mathbb{R}^{J_m} \to \mathbb{R}\cup\{\pm\infty\}$ satisfies that $g_m(\cdot, \bm{z}_m)$ is convex on $\mathbb{R}^N$ for any $\bm{z}_m\in\mathcal{Z}_m$ and $g_m(\bm{x}, \cdot)$ is concave on $\mathbb{R}^{J_m}$ for any $\bm{x}\in\mathcal{X}$, with $\mathcal{X}\times \mathcal{Z}_m \subseteq \mathrm{dom}(g_m)$. Here, $\mathrm{dom}(\cdot)$ denotes the domain.
		
\item\label{Assump: Robust-Prob-1-iii} {\normalfont (Existence of Slater Points)} There exists $\bar{\bm{x}} \in \mathcal{X}$ such that $\displaystyle \max_{\bm{z}_m \in \mathcal{Z}_m} g_m (\bar{\bm{x}}, \bm{z}_m) <0$ for any $m\in [M]$.
		
\item\label{Assump: Robust-Prob-1-iv} {\normalfont (Existence of Optimal Solutions)} The \ref{opt:RO1} problem has an optimal solution. 
\end{enumerate}
\end{assumption}

An immediate consequence of Assumption~\ref{Assump: Robust-Prob-1} is the following equivalence.

\begin{proposition}\label{prop:mmm-equivalence}
Suppose that Assumption~\ref{Assump: Robust-Prob-1} holds. 
Then the \ref{opt:RO1} problem is equivalent to the~\ref{opt:mmm1} problem
in the sense that their optimal values are equal and that for any optimal solution $(\bm{\lambda}^\star, \bm{x}^\star, \bm{z}^\star)$ to \ref{opt:mmm1}, $\bm{x}^\star$ is an optimal solution to \ref{opt:RO1}.
\end{proposition}

We also need the following assumption concerning the subgradient of the functions $g_1,\dots, g_m$ and $f_0$, which is customary in the literature of subgradient-type algorithms.

\begin{assumption}[Uniformly Bounded Subdifferentials]\label{Assump: Robust-Prob-2}
The function $f_0$ is subdifferentiable on $\mathcal{X}$. For any $m\in [M]$, $\bm{x}\in\mathcal{X}$ and $\bm{z}_m\in \mathcal{Z}_m$, the functions $g_m(\cdot, \bm{z}_m)$ and $-g_m(\bm{x}, \cdot)$ are subdifferentiable on $\mathcal{X}$ and $\mathcal{Z}_m$, respectively. 
There exist constants $D_0, D_1,\dots, D_M, E_1,\dots, E_M >0$ such that for any $m\in[M]$, $\bm{z}_m\in\mathcal{Z}_m$ and $\bm{x}\in\mathcal{X}$, 
\begin{equation*}
\begin{array}{r@{\;}l@{\;}l}
\|\bm{\xi}_0\|_2 \le & D_0 &\quad\forall \bm{\xi}_0\in \partial f_0(\bm{x}), \\
\|\bm{\xi}_m\|_2 \le & D_m &\quad\forall \bm{\xi}_m\in \partial_{\bm{x}} g_m(\bm{x}, \bm{z}_m), \\
\|\bm{\zeta}_m\|_2 \le & E_m &\quad\forall \bm{\zeta}_m\in \partial_{\bm{z}_m} (-g_m)(\bm{x}, \bm{z}_m).
\end{array}
\end{equation*}
\end{assumption}

Recall that a function is subdifferentiable at a point if its subdifferential (\ie, the set of subgradients) at that point is non-empty~\cite[Section~23]{rockafellar2015convex}. It is well-known that any convex function has a bounded non-empty subdifferential at any point in the interior of its domain~\cite[Theorem~23.4]{rockafellar2015convex}. Therefore, in view of Assumption~\ref{Assump: Robust-Prob-1}\ref{Assump: Robust-Prob-1-ii}, Assumption~\ref{Assump: Robust-Prob-2} can only be violated on the boundary and hence very mild. Indeed, if $\mathcal{X}\subseteq \mathrm{int}(\mathrm{dom}(f_0))$ and $\mathcal{X}\times \mathcal{Z}_m \subseteq \mathrm{int}(\mathrm{dom}(g_m))$ for all $m\in[M]$ (\eg, when $f_0$ and $g_1,\dots,g_m$ are real-valued everywhere), where $\mathrm{int}(\cdot)$ denotes the interior, then Assumption~\ref{Assump: Robust-Prob-2} holds.  

% what Assumption~\ref{Assump: Robust-Prob-2}\ref{Assump: Robust-Prob-2-i} additionally requires is the subdifferentiability of the functions on the boundary of the domains $\mathcal{Z}_1,\dots, \mathcal{Z}_M$ and $\mathcal{X}$. Furthermore, it is also well-known that the subdifferential is compact at any point in the relative interior~\cite[Theorem~23.4]{rockafellar2015convex}.
% Hence, Assumption~\ref{Assump: Robust-Prob-2}\ref{Assump: Robust-Prob-2-ii} can only be violated on the boundary of the sets $\mathcal{Z}_1,\dots, \mathcal{Z}_M$ and $\mathcal{X}$. For example, if $f_0$ is the indicator function on $\mathcal{X}$, then the subdifferential at any boundary point $\bm{x}$ is the normal cone to $\mathcal{X}$ at $\bm{x}$, which is unbounded. Thus, only minimal additional assumptions are imposed, as compared with the standard assumptions in the RO literature (\ie, Assumption~\ref{Assump: Robust-Prob-1}).

\subsection{Outer Convergence Rate}\label{subsec:out_rate}

Our first main theoretical result concerns the convergence rate of the outer algorithm.

\begin{theorem}\label{thm: iteration-complexity-Algo1}
Suppose that Assumptions~\ref{Assump: Robust-Prob-1} and~\ref{Assump: Robust-Prob-2} hold. Consider \Cref{alg:outer} with $\theta = \nu = \tfrac{1}{K}$, $\alpha \leq \frac{1}{\sqrt{\sum_{m\in [M]} D_m^2}}$ and $\beta \leq  \frac{1}{2\sqrt{\sum_{m\in [M]} D_m^2}}$. Then, the output $\bar{\bm{x}}^K$ satisfies 
\begin{equation*}
f_0(\bar{\bm{x}}^K) - f_0(\bm{x}^{\star}) \le \tfrac{C_{o,1}}{K} \quad\text{and}\quad \max_{m\in[M]} \left[ f_m(\bar{\bm{x}}^K) \right]_+ \le \tfrac{C_{o,2}}{K}
\end{equation*}
for some constants $C_{o,1}, C_{o,2}> 0$.
\end{theorem}

Although the outer algorithm shares a similar blueprint of the alternating proximal algorithm, \Cref{thm: iteration-complexity-Algo1} does not follow directly from existing studies but requires certain new ideas. First, even though alternating proximal algorithms with inaccurate updates have been studied in a number of works, the forms of inaccuracy assumed in those papers do not cover that of our outer algorithm. Second, existing theoretical works on alternating proximal algorithms primarily focus on bounding the saddle gap 
$\mathcal{L}( \bar{\bm{\lambda}}^K, \bm{x}^\star ) - \mathcal{L}( \bm{\lambda}^\star, \bar{\bm{x}}^K )$; see~\cite{chambolle2011first, chambolle2016ergodic}. We, however, care more about the optimality gap and constraint violation, since our ultimate goal is to solve the \ref{opt:RO1} problem.

\subsection{Inner Convergence Rate}\label{subsec:in_rate}

The following theorem asserts that under Assumptions~\ref{Assump: Robust-Prob-1} and~\ref{Assump: Robust-Prob-2}, the inner algorithm finds a strong $\nu$-approximate saddle point in $\mathcal{O}(\nu^{-2})$ iterations.

\begin{theorem}\label{thm:inner-rate}
Suppose that Assumptions~\ref{Assump: Robust-Prob-1} and~\ref{Assump: Robust-Prob-2} hold. There exists a constant $C_{i,1} > 0$ such that if $T \ge \frac{C_{i,1}}{\nu^2}$ and $\gamma , \delta \le \tfrac{1}{\sqrt{T}}$,
then the output $ ( \bar{\bm{x}}_T, \bar{\bm{z}}_T ) $ of \Cref{alg:inner} is a strong $\nu$-approximate saddle point.
\end{theorem}

To present our next result, we introduce the following smoothness conditions on the functions $f_0, g_1,\dots, g_M$.

\begin{assumption}\label{ass:smoothness}
The function $f_0$ is differentiable\footnote{A function is differentiable on a \textit{non-open} set $\mathcal{S}$ if it is differentiable on an open set containing $\mathcal{S}$.} on $\mathcal{X}$. For any $m\in[M]$, the function $g_m$ is differentiable on $\mathcal{X}\times \mathcal{Z}_m$. There exist constants $D'_0, D'_1,\dots, D'_M, E'_{1,1}, E'_{1,2},\dots, E'_{M,1}, E'_{M,2} >0$ such that for any $m\in[M]$, $\bm{z}_m, \bm{z}_m' \in \mathcal{Z}_m$ and $\bm{x},\bm{x}'\in\mathcal{X}$, we have
\begin{equation*}
\begin{array}{r@{\;}l}
    \| \nabla f_0( \bm{x} ) - \nabla f_0(\bm{x}' ) \|_2 \le & D_0' \, \| \bm{x} - \bm{x}'  \|_2 , \\
    \| \nabla_{\bm{x}} g_m( \bm{x}, \bm{z}_m ) - \nabla_{\bm{x}} g_m(\bm{x}', \bm{z}_m ) \|_2 \le & D'_m \, \| \bm{x} - \bm{x}' \|_2 , \\
    \| \nabla_{\bm{z}_m} g_m( \bm{x}, \bm{z}_m ) - \nabla_{\bm{z}_m} g_m(\bm{x}', \bm{z}_m' ) \|_2 \le & E'_{m,1} \, \| \bm{x} - \bm{x}' \|_2 + E'_{m,2} \, \| \bm{z}_m - \bm{z}_m' \|_2 .
\end{array}
\end{equation*}
\end{assumption}

Under Assumption~\ref{Assump: Robust-Prob-1}, it can be readily shown that Assumption~\ref{ass:smoothness} implies Assumption~\ref{Assump: Robust-Prob-2}.
The theorem below asserts that the iteration complexity of the inner algorithm can be improved to $\mathcal{O}(\nu^{-1})$ if we replace Assumption~\ref{Assump: Robust-Prob-2} by Assumption~\ref{ass:smoothness}.

\begin{theorem}\label{thm:inner-rate-smooth}
Suppose that Assumptions~\ref{Assump: Robust-Prob-1} and~\ref{ass:smoothness} hold. There exists a constant $C_{i,2} > 0$ such that if $T \ge  \frac{C_{i,2}}{\nu}$,
\begin{equation*}
    \gamma \le \frac{1}{ \displaystyle D_0'  + \sum_{m\in [M]} \lambda_m^{k+1} D'_m + \sqrt{2 \sum_{m\in[M]} \left(\lambda_m^{k+1} E'_{m,1}\right)^2 }}, 
\end{equation*}
and
\begin{equation*}
    \delta \le \frac{1}{ \displaystyle 2\sqrt{2} \left( \max_{m\in[M]} \lambda_m^{k+1} E'_{m,2}\right) + \sqrt{2 \sum_{m\in[M]} \left(\lambda_m^{k+1} E'_{m,1}\right)^2 }} ,
\end{equation*}
then the output $(\bar{\bm{x}}_T, \bar{\bm{z}}_T ) $ of \Cref{alg:inner} is a strong $\nu$-approximate saddle point.
\end{theorem}

\subsection{Oracle Complexity}

Since our approach assumes access to the subgradient oracles of the functions $f_0, g_1,\dots,g_M$ as well as the projection oracles of the sets $\mathcal{X}, \mathcal{Z}_1,\dots, \mathcal{Z}_M$, a faithful and popular measure of computational efficiency would be the so-called \textit{oracle complexity}---the total number of calls of these oracles---for achieving an $\varepsilon$-approximate optimal solution. Here we recall that for any $\varepsilon > 0$, a point $\bm{x}$ is an $\varepsilon$-approximate optimal solution to \ref{opt:RO1} if
\[ 
f_0(\bm{x} ) - f_0 (\bm{x}^\star) \le \varepsilon \quad\text{and}\quad \max_{m\in [M]} \; [f_m (\bm{x})]_+ \le \varepsilon. 
\]
The next theorem presents the oracle complexity of \prom by combining the convergence rate results for the outer and inner algorithms.

\begin{theorem}\label{thm:complexity_1}
Suppose that Assumptions~\ref{Assump: Robust-Prob-1} and~\ref{Assump: Robust-Prob-2} hold. Then, for any $\varepsilon>0$, the oracle complexity of \prom for achieving an $\varepsilon$-approximate optimal solution to \ref{opt:RO1} is $\mathcal{O}(\varepsilon^{-3})$.
\end{theorem}

Under similar assumptions as in Theorem~\ref{thm:complexity_1}, \cite{postek2021first} proved that the oracle complexity of the algorithm SGSP is $\mathcal{O}(\varepsilon^{-2})$. However, due to the use of perspective transformation, SGSP relies on subgradient and projection oracles that are generally more computationally expensive. Thus, our oracle complexity results are not directly comparable to that of \cite{postek2021first}.

The oracle complexity can be improved if Assumption~\ref{Assump: Robust-Prob-2} is replaced by Assumption~\ref{ass:smoothness}.

\begin{theorem}\label{thm:complexity_2}
Suppose that Assumptions~\ref{Assump: Robust-Prob-1} and~\ref{ass:smoothness} hold. Then, for any $\varepsilon>0$, the oracle complexity of \prom for achieving an $\varepsilon$-approximate optimal solution to \ref{opt:RO1} is $\mathcal{O}(\varepsilon^{-2})$.
\end{theorem}

Under similar assumptions to Theorem~\ref{thm:complexity_2}, a first-order method is developed in \cite{ho2018online, ho2019exploiting} via online convex optimization and proved to enjoy the oracle complexity $\mathcal{O}(\varepsilon^{-1}\log \tfrac{1}{\varepsilon})$. Nevertheless, the authors considered only low- to medium-dimensional instances in their numerical experiments. In Section~\ref{sec:exp}, we demonstrate that when compared against the first-order methods in \cite{postek2021first} and \cite{ho2018online, ho2019exploiting} our algorithm \prom is substantially more stable and efficient.

\section{Extension to Projection-Unfriendly Uncertainty Sets}\label{sec:extension}

For some applications of RO, the uncertainty sets $\mathcal{Z}_m$ take the form of an intersection and do not admit an easy projection. Directly invoking our \prom in Section~\ref{sec:prom} to such RO problems could be inefficient. Below we extend our first-order algorithm \prom to RO problems with uncertainty sets of the form
\begin{equation}\label{def:generalized-Zm}
\mathcal{Z}_m = \widetilde{\mathcal{Z}}_m \cap \left(\bigcap_{i\in [I_m]} \mathcal{Z}_{m,i} \right),
\end{equation}
where for each $i \in [I_m]$, the set $\mathcal{Z}_{m,i}  =  \{\bm{z}_m \in \mathbb{R}^{J_m}\mid h_{m, i}(\bm{z}_m)\leq  0\}$ is defined by some function $h_{m,i}$.

To present the extended \prom and its analysis, we make the following assumption: an adaptation of Assumption~\ref{Assump: Robust-Prob-1} to the projection-unfriendly setting.

\begin{assumption}\label{Assump: Robust-Prob-1-gen}
The following conditions hold.
\begin{enumerate}[font=\normalfont, label=(\textit{\roman*})]
\item\label{Assump: Robust-Prob-1-gen-i} {\normalfont (Compactness and Convexity of Sets)} The sets $\widetilde{\mathcal{Z}}_1\subseteq \mathbb{R}^{J_1},\dots, \widetilde{\mathcal{Z}}_M\subseteq \mathbb{R}^{J_M}$ and $\mathcal{X}\subseteq \mathbb{R}^N$ are non-empty, compact and convex.
		
\item\label{Assump: Robust-Prob-1-gen-ii} {\normalfont (Convexity of Functions)} 
The function $f_0: \mathbb{R}^N \to \mathbb{R}\cup\{+\infty\}$ is convex, with $\mathcal{X}\subseteq \mathrm{dom}(f_0)$. 
For any $m\in[M]$, the function $g_m:\mathbb{R}^N\times \mathbb{R}^{J_m} \to \mathbb{R}\cup\{\pm\infty\}$ satisfies that $g_m(\cdot, \bm{z}_m)$ is convex on $\mathbb{R}^N$ for any $\bm{z}_m\in\widetilde{\mathcal{Z}}_m$ and $g_m(\bm{x}, \cdot)$ is concave on $\mathbb{R}^{J_m}$ for any $\bm{x}\in\mathcal{X}$, with $\mathcal{X}\times \widetilde{\mathcal{Z}}_m \subseteq \mathrm{dom}(g_m)$. 
For any $m\in [M]$ and $i\in [I_m]$, the function $h_{m,i}: \mathbb{R}^{J_m} \to \mathbb{R}\cup\{+\infty\}$ is convex, with $\widetilde{\mathcal{Z}}_m\subseteq \mathrm{dom} (h_{m,i})$.
		
\item\label{Assump: Robust-Prob-1-gen-iii} {\normalfont (Existence of Slater Points)} There exists $\bar{\bm{x}} \in \mathcal{X}$ such that $\displaystyle \max_{\bm{z}_m \in \widetilde{\mathcal{Z}}_m} g_m (\bar{\bm{x}}, \bm{z}_m) <0$ for any $m\in[M]$. For any $m\in[M]$, there exists $\bar{\bm{z}}_m\in \widetilde{\mathcal{Z}}_m$ such that $\displaystyle h_{m,i} (\bar{\bm{z}}_m) <0$ for any $i\in [I_m]$.
		
\item\label{Assump: Robust-Prob-1-gen-iv} {\normalfont (Existence of Optimal Solutions)} The \ref{opt:RO1} problem has an optimal solution. 
\end{enumerate}
\end{assumption}

We also need an adaptation of Assumption~\ref{Assump: Robust-Prob-2}.
\begin{assumption}[Uniformly Bounded Subdifferentials]\label{Assump: Robust-Prob-2-gen}

The function $f_0$ is subdifferentiable on $\mathcal{X}$. For any $m\in [M]$, $\bm{x}\in\mathcal{X}$ and $\bm{z}_m\in \widetilde{\mathcal{Z}}_m$, the functions $g_m(\cdot, \bm{z}_m)$ and $-g_m(\bm{x}, \cdot)$ are subdifferentiable on $\mathcal{X}$ and $\widetilde{\mathcal{Z}}_m$, respectively. For any $m\in [M]$ and $i\in [I_m]$, the function $h_{m,i}$ is subdifferentiable on $\widetilde{\mathcal{Z}}_m$.
There exist constants $D_0, D_1,\dots, D_M, E_1,\dots, E_M, F_1,\dots, F_M >0$ such that for any $m\in[M]$, $ i\in [I_m]$, $\bm{z}_m\in\widetilde{\mathcal{Z}}_m$ and $\bm{x}\in\mathcal{X}$, we have
\begin{equation*}
\begin{array}{r@{\;}l@{\;}l}
    \|\bm{\xi}_0\|_2 \le & D_0  &\quad\forall \bm{\xi}_0\in \partial f_0(\bm{x}),  \\
    \|\bm{\xi}_m\|_2 \le & D_m &\quad\forall \bm{\xi}_m\in \partial_{\bm{x}} g_m(\bm{x}, \bm{z}_m) , \\
    \|\bm{\zeta}_m\|_2 \le & E_m & \quad\forall \bm{\zeta}_m\in \partial_{\bm{z}_m} (-g_m)(\bm{x}, \bm{z}_m) , \\
    \|\bm{\eta}_{m,i}\|_2 \le & F_m &\quad\forall \bm{\eta}_{m,i} \in \partial h_{m,i}(\bm{z}_m).
\end{array}
\end{equation*}
\end{assumption}

Similarly, the oracle complexity of the extended \prom can be improved when the gradients of the constituent functions satisfy certain Lipschitz property.
\begin{assumption}\label{ass:smoothness-gen}
The function $f_0$ is differentiable on $\mathcal{X}$. For any $m\in[M]$, the function $g_m$ is differentiable on $\mathcal{X}\times\widetilde{\mathcal{Z}}_m$. For any $m\in[M]$ and $i\in [I_m]$, the function $h_{m,i}$ is differentiable on $\widetilde{\mathcal{Z}}_m$. There exist constants $D'_0, D'_1,\dots, D'_M, E'_{1,1}, E'_{1,2},\dots, E'_{M,1}, E'_{M,2}, F'_1,\dots, F'_M >0$ such that for any $m\in[M]$, $ i\in [I_m]$, $\bm{z}_m, \bm{z}_m' \in\widetilde{\mathcal{Z}}_m$ and $\bm{x},\bm{x}'\in\mathcal{X}$, we have
\begin{equation*}
\begin{array}{r@{\;}l}
    \| \nabla f_0( \bm{x} ) - \nabla f_0(\bm{x}' ) \|_2 \le & D_0' \, \| \bm{x} - \bm{x}'  \|_2 , \\
    \| \nabla_{\bm{x}} g_m( \bm{x}, \bm{z}_m ) - \nabla_{\bm{x}} g_m(\bm{x}', \bm{z}_m ) \|_2 \le & D'_m \, \| \bm{x} - \bm{x}' \|_2 , \\
    \| \nabla_{\bm{z}_m} g_m( \bm{x}, \bm{z}_m ) - \nabla_{\bm{z}_m} g_m(\bm{x}', \bm{z}_m' ) \|_2 \le & E'_{m,1} \, \| \bm{x} - \bm{x}' \|_2 + E'_{m,2} \, \| \bm{z}_m - \bm{z}_m' \|_2, \\
    \| \nabla h_{m,i}( \bm{z}_m ) - \nabla h_{m,i}(\bm{z}_m' ) \|_2 \le & F_m' \, \| \bm{z}_m - \bm{z}_m'  \|_2 .
\end{array}
\end{equation*}
\end{assumption}

Note that for any $m\in[M]$, $g_m (\cdot, \bar{\bm{z}}_m)$ is continuous on $\mathcal{X}$ because it is real-valued and convex on $\mathcal{X}$. By the compactness of $\mathcal{X}$, there exist constants $G_1,\dots, G_M < 0 $ such that
\begin{equation*}
     g_m ( \bm{x}, \bar{\bm{z}}_m ) \ge G_m \qquad\forall \bm{x}\in\mathcal{X}, \; m\in[M].
\end{equation*}
For simplicity, we denote $\bm{h}_m (\bm{z}_m) = (h_{m,1}(\bm{z}_m), \dots, h_{m,I_m}(\bm{z}_m) ) $, $\widetilde{\mathcal{Z}} = \widetilde{\mathcal{Z}}_1\times\cdots\times\widetilde{\mathcal{Z}}_M$ and $\mathcal{M} = [0,a_1]^{I_1} \times\cdots\times [0,a_M]^{I_M}$, 
where
\[ 
a_m = \frac{G_m}{\max_{i\in [I_{m}]} \{ h_{m, i}(\bar{\bm{z}}_m)\}}.  
\]

Our extension of \prom to projection-unfriendly uncertainty sets of the form~\eqref{def:generalized-Zm} is based on the following proposition.
\begin{proposition}\label{prop:extended_RO}
    Consider the \ref{opt:RO1} problem with uncertainty sets $\mathcal{Z}_m$ as in~\eqref{def:generalized-Zm}. Suppose that Assumption~\ref{Assump: Robust-Prob-1-gen} holds. Then, the \ref{opt:RO1} problem is equivalent to the problem
    \begin{equation}\label{opt:Gene-RO1}\tag{\sc $\widetilde{\text{Robust}}$}
		\begin{array}{rll}
				\min & \tilde{f}_0(\bmt{x}) \\
				{\rm s.t.} & \displaystyle  \max_{\bm{z}_m \in \widetilde{\mathcal{Z}}_m} \tilde{g}_m (\bmt{x}, \bm{z}_m) \leq 0  \quad\forall m\in [M] \\
				& \bmt{x}=(\bm{x}, \bm{\mu}_1,\dots,\bm{\mu}_M ) \in \widetilde{\mathcal{X}},
			\end{array}
	\end{equation}
where $ \tilde{f}_0(\bmt{x}) = f_0 (\bm{x}) $, $\widetilde{\mathcal{X}} = \mathcal{X}\times \mathcal{M}$ and $\tilde{g}_m(\bmt{x}, \bm{z}_m) = g_m (\bm{x}, \bm{z}_m)- \bm{\mu}_m^\top \bm{h}_m (\bm{z}_m) $ for all $m\in [M]$.
\end{proposition}

\begin{algorithm}[t]
    \caption{Extended Outer Algorithm. \label{alg:outer-extension}}
    \SetKwInOut{Input}{Input}
    \SetKwInOut{Output}{Output}
    \Input{$K\geq 1$, $\beta>0$, $\alpha>0$, $\nu>0$, $\theta>0$, $\bm{\lambda}^{0} = \bm{0}$, $ \bm{x}^{0}\in \mathcal{X}$ and $\bm{\mu}^0\in \mathcal{M}$.}
    \For{$k= 0,1, \dots, K-1$}{
        ($\bm{\lambda}$-update) For all $m\in [M]$, find $\bm{z}_{m}^{k} \in \widetilde{\mathcal{Z}}_m$ satisfying
        \[
            \max_{\bm{z}_m \in \widetilde{\mathcal{Z}}_m}  g_m (\bm{x}^k, \bm{z}_m)-(\bm{\mu}_m^{k})^{\top} \bm{h}_m(\bm{z}) \le g_{m}(\bm{x}^{k},~ \bm{z}_{m}^{k}) - (\bm{\mu}_m^k)^{\top} \bm{h}_m(\bm{z}_m^k) + \theta,
        \]
        and set
            \begin{equation*}
            \lambda_m^{k+1} = \left[\lambda_m^k + \beta \left(2\, g_m(\bm{x}^k, \bm{z}_m^k)-2(\bm{\mu}_m^{k})^{\top} \bm{h}_m(\bm{z}_m^k) - g_m(\bm{x}^{k-1}, \bm{z}_m^{k-1}) + (\bm{\mu}_m^{k-1})^{\top} \bm{h}_m(\bm{z}_m^{k-1})\right)  \right]_+ ,
            \end{equation*}
        where $\bm{x}^{-1} = \bm{x}^0$ and $\bm{z}^{-1} = \bm{z}^0$ for the $0$-th iteration. \vspace{3mm}
      
        ($\bmt{x}$-update)  Compute a strong $\nu$-approximate saddle point  $((\bm{x}^{k+1}, \bm{\mu}^{k+1}),  \bmt{z}^{k+1})$ of the following problem via the Extended Inner Algorithm (Algorithm~\ref{alg:inner-extension}):
        \begin{equation}\label{def:G-inner}\tag{\sc $\widetilde{\text{Inner}}$}
        \begin{split}
            \min_{\bm{x} \in \mathcal{X},\, \bm{\mu} \in \mathcal{M}} \; \max_{\bm{z} \in \widetilde{\mathcal{Z}}} \; \Bigg\{ & f_0(\bm{x})+ \sum_{m\in [M] }\lambda_{m}^{k+1}  \left(  g_m(\bm{x}, \bm{z}_m)  - \bm{\mu}_m^{\top} \bm{h}_m(\bm{z}_m)\right) \\
            & \quad+ \frac{1}{2\alpha} \|\bm{x}- \bm{x}^k\|_2^2+\frac{1}{2\alpha}\|\bm{\mu} -\bm{\mu}^k\|_2^2 \Bigg\},
        \end{split}
        \end{equation}
        }
    \Output{$\bar{\bm{x}}^{K} =\frac{1}{K} \sum_{k\in [K] }\bm{x}^{k}$.}
\end{algorithm}

The~\ref{opt:Gene-RO1} problem is obtained by penalizing the constraints $\bm{h}_m (\bm{z}_m) \le \bm{0}$ in the $m$-th embedded problem to its objective. It is expected to be equivalent to the original~\ref{opt:RO1} problem under suitable assumptions if we do not restrict the dual variable $\bm{\mu}_m$ from above, \ie, if $a_m $ in the definition of $\mathcal{M}$ is replaced by $+\infty$ for all $m\in [M]$. Nevertheless, for our algorithmic framework and theoretical results in Section~\ref{sec:prom} to be applicable, we need to introduce an upper bound $a_m$ for each $\bm{\mu}_m$ to make the feasible region $\widetilde{\mathcal{X}}$ compact. What is perhaps less trivial is that the equivalence remains after restricting the dual variable.

Applying our framework to the \ref{opt:Gene-RO1} problem, we obtain an extension of \prom for RO problems with projection-unfriendly uncertainty sets of the form~\eqref{def:generalized-Zm}. To facilitate easy usage, we present the extended outer and inner algorithms fully in terms of the basic constituent functions and sets in Algorithms~\ref{alg:outer-extension} and~\ref{alg:inner-extension}, respectively. The extended \prom enjoys the following complexity result.
 
\begin{algorithm}[t]
	\caption{Extended Inner Algorithm. \label{alg:inner-extension}}
	\SetKwInOut{Input}{Input}
	\SetKwInOut{Output}{Output}
	\Input{$T\ge 1$,  $\delta >0$,  $\gamma>0$, $\bm{x}_0, \bm{x}^k\in \mathcal{X}$, $ \bm{\mu}_0, \bm{\mu}^k\in \mathcal{M}$, $\bm{z}_0\in \widetilde{\mathcal{Z}}$ and $\bm{\lambda}^{k+1} \ge \bm{0}$.}
	\For{$t = 0, \dots, T-1$}
	{
		($\bm{z}$-update) Compute $\bm{\zeta}_{t,m} \in \partial_{\bm{z}_m} (-g_m)( \bm{x}_t, \bm{z}_{t,m} )$ and $\bm{\eta}_{t,m,i}\in\partial h_{m,i}( \bm{z}_{t,m} )$ for all $m\in [M]$ and $ i\in [I_m]$. Set
		\[ 
		\bm{z}_{t+1,m} = \Proj_{\widetilde{\mathcal{Z}}_m} \left( \bm{z}_{t,m} - \delta \,\lambda_m^{k+1}\left( 2\bm{\zeta}_{t,m} - \bm{\zeta}_{t-1,m} + \sum_{i = 1}^{I_m} \mu_{m,i} (2 \bm{\eta}_{t,m,i} - \bm{\eta}_{t-1,m,i}) \right) \right), 
		\]
		where $\bm{\zeta}_{-1,m} = \bm{\zeta}_{0,m}$ and $\bm{\eta}_{-1,m,i} = \bm{\eta}_{0,m,i}$ for the $0$-th iteration. \vspace{3mm}
		
		($\bmt{x}$-update) Compute $\bm{\xi}_{t,0} \in \partial f_0 (\bm{x}_t)$ and $\bm{\xi}_{t,m} \in \partial_{\bm{x}} g (\bm{x}_t, \bm{z}_{t+1, m})$ for all $m\in[M]$. Set $\bm{\xi}_t = \bm{\xi}_{t,0} + \lambda_1^{k+1}\, \bm{\xi}_{t,1}+\cdots+ \lambda_M^{k+1}\, \bm{\xi}_{t,M}$,
		\begin{equation*}
			    \begin{split}
                \bm{x}_{t+1}  =&  \Proj_{\mathcal{X}} \left( \frac{\alpha \gamma}{\alpha + \gamma}\left(\frac{1}{\alpha}\bm{x}^k + \frac{1}{\gamma} \bm{x}_t - \bm{\xi}_t\right) \right) \quad {\rm and} \quad \\
                \bm{\mu}_{t+1,m}  =&  \Proj_{[0, a_m]^{I_m}} \left( \frac{\alpha \gamma}{\alpha + \gamma}\left(\frac{1}{\alpha}\bm{\mu}_m^k + \frac{1}{\gamma} \bm{\mu}_{t,m} + \lambda_m^{k+1} \bm{h}_m(\bm{z}_{t+1,m})\right) \right).
				\end{split}
		\end{equation*}
		}
	\Output{$ (\bar{\bm{x}}_T, \bar{\bm{\mu}}_T) = (\frac{1}{T}\sum_{t\in [T]}\bm{x}_{t}, \frac{1}{T}\sum_{t\in [T]}\bm{\mu}_{t}) $ and $\bar{\bm{z} }_T = \frac{1}{T}\sum_{t\in [T]}\bm{z}_{t}$.}
\end{algorithm}

\begin{theorem}
    \label{thm:extended_prom_complexity}
    Consider the \ref{opt:RO1} problem with uncertainty sets $\mathcal{Z}_m$ as in~\eqref{def:generalized-Zm}. Suppose that Assumptions~\ref{Assump: Robust-Prob-1-gen} and~\ref{Assump: Robust-Prob-2-gen} hold.
    Then, for any $\varepsilon >0$, the oracle complexity of extended \prom (Algorithms~\ref{alg:outer-extension} and~\ref{alg:inner-extension} combined) for achieving an $\varepsilon$-approximate optimal solution to \ref{opt:RO1} is $\mathcal{O}(\varepsilon^{-3})$. If Assumption~\ref{Assump: Robust-Prob-2-gen} is replaced by Assumption~\ref{ass:smoothness-gen}, then the oracle complexity is strengthened to $\mathcal{O}(\varepsilon^{-2})$. 
\end{theorem}

A major part of the proof of Theorem~\ref{thm:extended_prom_complexity} is to verify that problem~\ref{opt:Gene-RO1} satisfies all the assumptions required by our results in Section~\ref{sec:convergence_analysis} (\ie, Assumptions~\ref{Assump: Robust-Prob-1}-\ref{ass:smoothness}).

\section{Numerical Experiments}\label{sec:exp}

This section explores the practical performance of the proposed algorithm \prom through extensive numerical experiments. We compare its performance with the reformulation approach, the cutting-plane method, and two recently developed first-order methods by~\cite{ho2018online} and~\cite{postek2021first}. We use the legends ``CP'' for the cutting-plane method, ``OCO'' for the first-order method in~\cite{ho2018online}, and ``SGSP'' for the first-order method in~\cite{postek2021first}. For the reformulation approach, we use Ref-$\varepsilon$, where $\varepsilon$ is the stopping accuracy and set to the usual accuracy level of first-order methods, $ 10^{-4}$ or $10^{-5}$, for a fair comparison. Furthermore, when we calculate the time for the reformulation approach, we calculate only the solver time but exclude the modeling and compilation time due to the interfacing. All algorithms are implemented in Python, and all experiments are performed on a $2.6$GHz laptop with $16$GB memory. 

\subsection{Robust QCQP}\label{sec:robust_QCQP}
Our first experiment concerns the following robust quadratically constrained quadratic program appeared in \cite{ho2018online,postek2021first}:
\begin{equation}\label{prob:RO-ex1-1}
\begin{array}{rll}
\min & \displaystyle \max_{\bm{z}_0\in \mathcal{Z}_0} g_0(\bm{x}, \bm{z}_0)\\
{\rm s.t.} & \displaystyle \max_{\bm{z}_m \in \mathcal{Z}_m} g_m(\bm{x}, \bm{z}_m) \leq 0 &~\forall m \in [M] \\
& \bm{x} \in \mathcal{X},
\end{array}
\end{equation}
where $\mathcal{X} =\{ \bm{x}\in \mathbb{R}^{N}\mid \|\bm{x}\|_2 \leq 1\}$ and for $m\in \{0\}\cup [M]$, $\mathcal{Z}_{m} = \{\bm{z}_m \in \mathbb{R}^{J_m} \mid \|\bm{z}_m\|_2 \leq 1\}$,
\[
\begin{split}
g_{m}(\bm{x}, \bm{z}_m)=\left\|\left(\bm{P}_{m 0}+\sum_{j=1}^{J_{m}} \bm{P}_{mj} z_{mj}\right) \bm{x}\right\|_2^{2}+  \bm{b}_{m}^{\top}\bm{x} +c_{m}.
\end{split}
\]
Here, $\bm{b}_m\in \mathbb{R}^{N}$, $c_m\in \mathbb{R}$, $\bm{P}_{mj}\in \mathbb{R}^{P \times N}$ and $z_{mj}$ is the $j$-th entry of $\bm{z}_m$ for all $j \in [J_m]$.
We generate the problem data $\bm{P}_{m j}$, $\bm{b}_{m}$ and $c_m$ in the same manner as in \cite{postek2021first}. For all $m$ and $j$, the entries of $\bm{P}_{mj}$ and $\bm{b}_{m}$ are i.i.d. uniform random variables on $[-1,1]$, normalized via $\bm{P}_{mj} \leftarrow \bm{P}_{m j}/\|[\bm{P}_{m0}^{\top} \cdots \bm{P}_{mK}^{\top} ]^{\top}\|_{2}$ and $\bm{b}_m \leftarrow \bm{b}_m/\|\bm{b}_{m}\|_2$. We fix $c_m=-0.05$ to ensure that problem~\eqref{prob:RO-ex1-1} has a Slater point.
Note that the each $g_m (\bm{x}, \bm{z}_m)$ is convex in $\bm{x}$ but not concave in $\bm{z}_m$. Nonetheless, by using the techniques in \cite{ho2018online,postek2021first}, we can transform problem~\eqref{prob:RO-ex1-1} into an instance of \ref{opt:RO1} satisfying all required assumptions. Note also that the reformulation in this case is a semidefinite program~\cite[Theorem~3.2]{ben1998robust}.

\begin{figure}[tb]
\footnotesize
\centering
\caption{$(M,N,P,J_{m}) = (3,1500,30,30)$.}
\label{fig-Medium-case-M=3_N=1500_I=30}
\includegraphics[scale=0.5]
{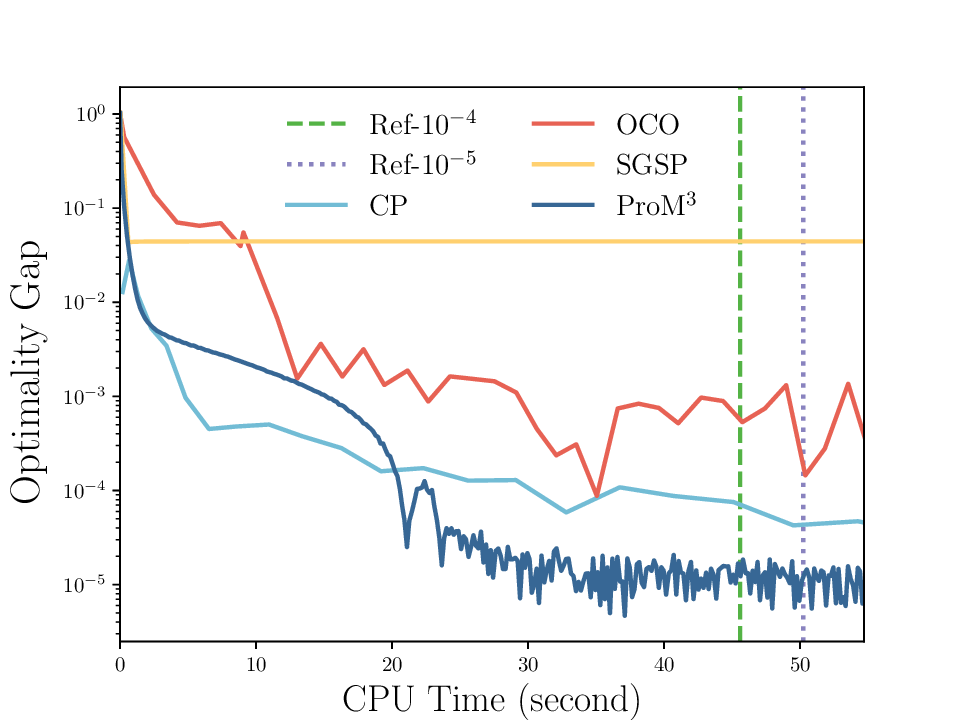}
\includegraphics[scale=0.5]{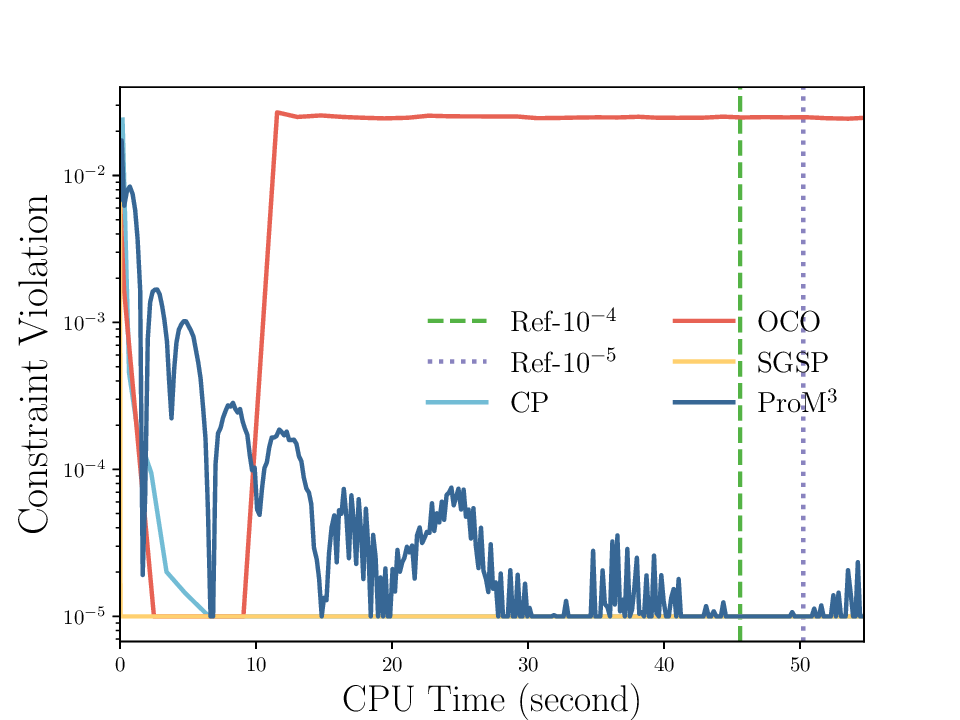}
\end{figure}

\begin{figure}[tb]
\footnotesize
\centering
\caption{{$(M,N,P,J_{m}) = (40,1500,30,30)$.}\label{fig-Medium-case-M=40_N=1500_I=30}}
\includegraphics[scale=0.5]{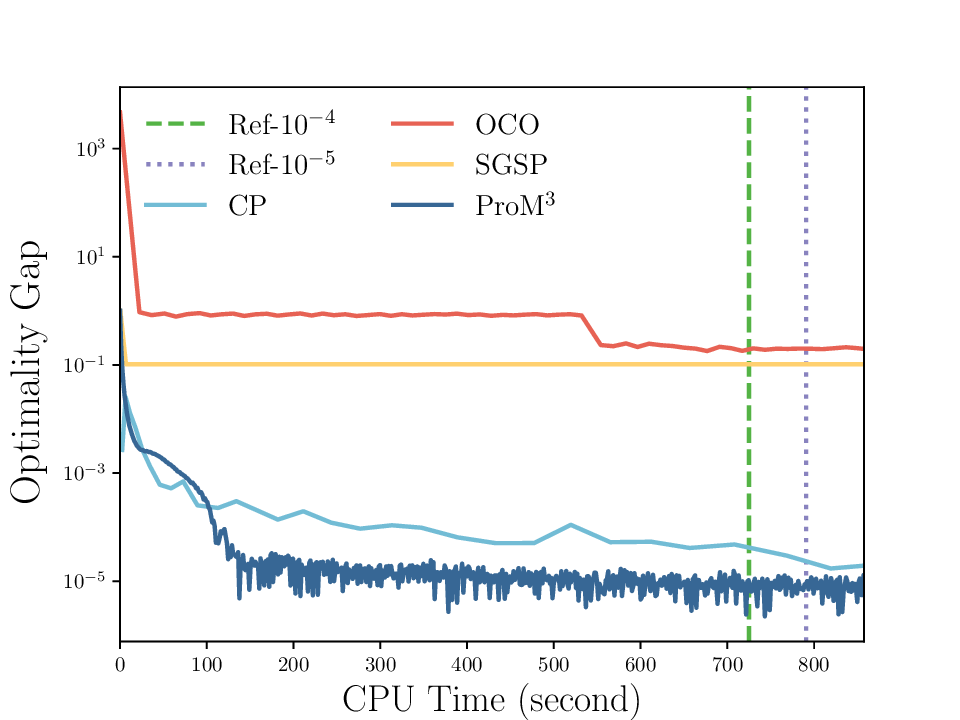}
\includegraphics[scale=0.5]{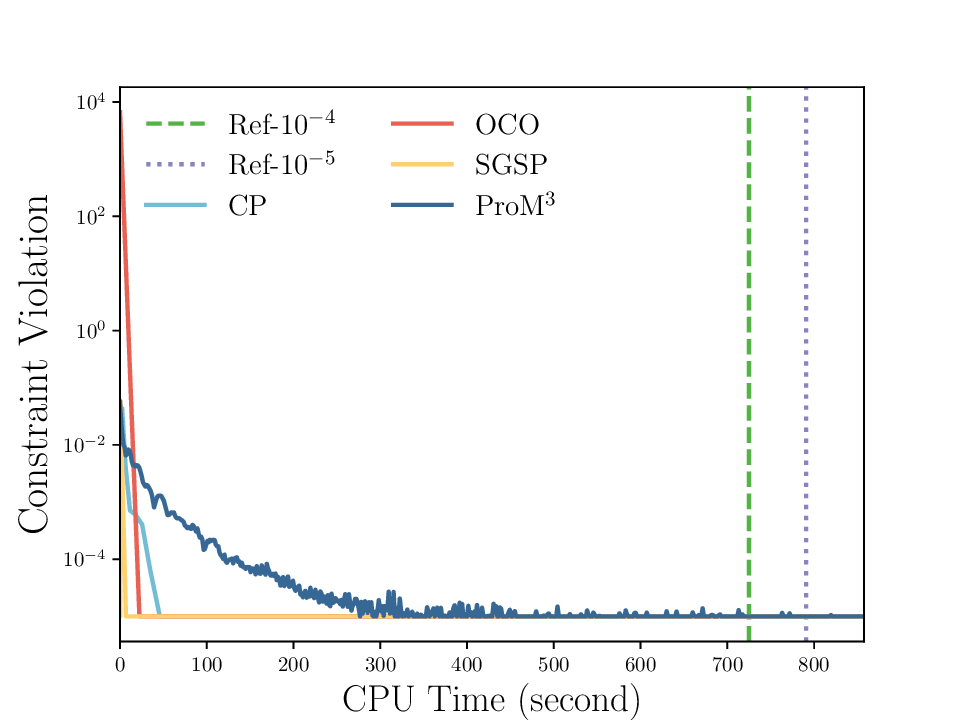}
\end{figure}

\begin{figure}[tb]
\footnotesize
\centering
\caption{{$(M,N,P,J_{m}) = (40,8000,50,50)$.}\label{fig-Medium-case-M=40_N=8000_I=50}}
\includegraphics[scale=0.5]
{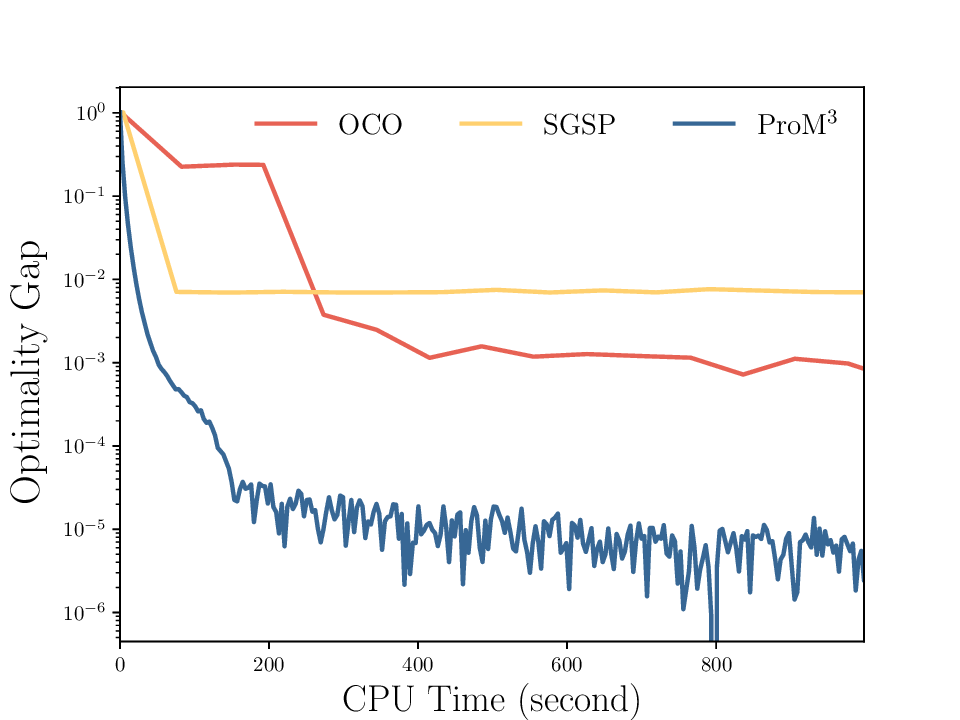}
\includegraphics[scale=0.5]
{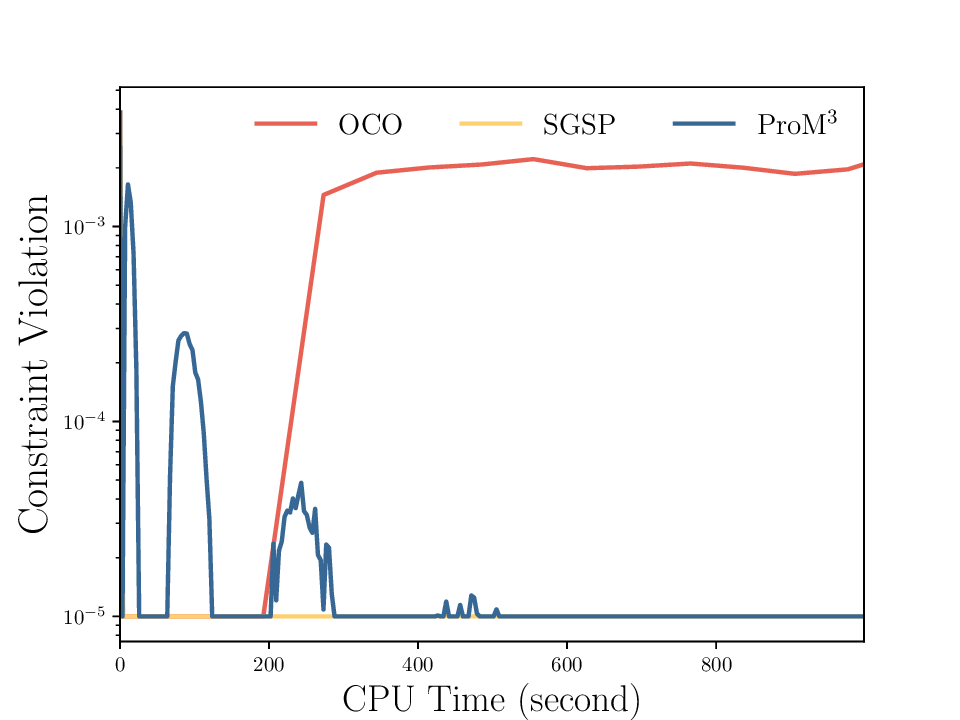}
\end{figure}

We test the algorithms on three problem dimensions: $(M,N,P,J_{m}) = (3,1500,30,30)$, $(M,N,P,J_{m}) = (40,1500,30,30)$ and $(M,N,P,J_{m}) = (4,8000,50,50)$, and the corresponding results are plotted in Figures~\ref{fig-Medium-case-M=3_N=1500_I=30}-\ref{fig-Medium-case-M=40_N=8000_I=50}, with the optimality gap $|f_0(\bm{x}) - f_0(\bm{x}^\star)|$ and the constraint violation $\max_{m\in [M]} [f_m (\bm{x})]_+$ shown on the left and right panels, respectively. 
To compute the optimality gap, for the first two cases, we take the solution returned by the reformulation approach with the default accuracy $10^{-12}$ as the ``true" optimum $\bm{x}^\star$. However, for the last case, the reformulation and cutting-plane methods do not work due to memory issues. In this case, we take the solution returned by our algorithm \prom as $\bm{x}^\star$, as it achieves a much higher accuracy than the other two first-order methods. Since we cannot keep track of the iterations of the solver in the reformulation approach, we only indicate its total time as a vertical line.

From Figures~\ref{fig-Medium-case-M=3_N=1500_I=30} and~\ref{fig-Medium-case-M=40_N=1500_I=30}, our algorithm \prom achieves the optimality gap $10^{-4}$ and $10^{-5}$ much faster than the reformulation approach and cutting-plane method, while the other two competing first-order methods get stuck at the level of $10^{-1}$ to $10^{-2}$. In terms of constraint violation, all the tested algorithms reach feasibility in a reasonable amount of time.
Figure~\ref{fig-Medium-case-M=40_N=8000_I=50} shows the result for the highest dimensional and the most challenging case, which is beyond the reach of the reformulation approach and cutting-plane method. Our algorithm \prom again considerably outperforms the two other first-order methods.

\subsection{Robust Log-Sum-Exponential Constraint}\label{sec:log-sum-exp}

We then consider a more involved RO problem with a highly non-linear embedded optimization problem. This should be contrasted with the robust QCQP example in Section~\ref{sec:robust_QCQP}, where the embedded problem's objective function and feasible region (uncertainty set) are both quadratic. Specifically, we consider the problem
\begin{equation}\label{problem-logistic}
\begin{array}{rll}
\min & \bm{c}^{\top}\bm{x} \\
{\rm s.t.} & \displaystyle \max_{\bm{z}_m \in \mathcal{Z}_m} g_m (\bm{x}, \bm{z}_m) \leq 0  \quad\forall m\in [M] \\
& \bm{x} \in \mathcal{X},
\end{array}
\end{equation}
where $\mathcal{X}=\{\bm{x}\mid -\bm{1} \leq \bm{x} \leq \bm{1}\}$ and for each $m\in[M]$, $\mathcal{Z}_m =\{\bm{z}_m \in \mathbb{R}^{J_m} \mid \bm{l} \leq \bm{z}_m \leq \bm{u}\}$ and
\[
g_m (\bm{x}, \bm{z}_m) = \bm{x}^{\top} \bm{A}_m \bm{z}_m  - d_m + \log\left(\bm{z}_{m,1} + \sum_{j=2}^{J_m} \bm{z}_{m,j} \exp{(\bm{b}_{m,j}^{\top} \bm{x}})\right).
\]
Such functions $g_m$ are called log-sum-exponential functions, and RO problems with robust log-sum-exponential constraints have been considered in \cite[Example 28]{ben2015deriving} and \cite{bertsimas2022robust}. The RO problem~\eqref{problem-logistic} is computationally very challenging since the robust log-sum-exponential constraint exhibits high non-linearity.
To the best of our knowledge, there is no tractable convex reformulation for the RO problem~\eqref{problem-logistic}~\cite[Example~28]{ben2015deriving}. We therefore do not compare with the reformulation approach in this experiment. Nevertheless, all other methods can be applied. This also shows that the scope of the reformulation approach is in general more restricted.
\begin{figure}[tb]
\centering
\caption{$(M, N , J_m) = (5, 200, 1000)$. \label{fig-ex2-N=200_q=1000}}
\includegraphics[scale=0.5]
{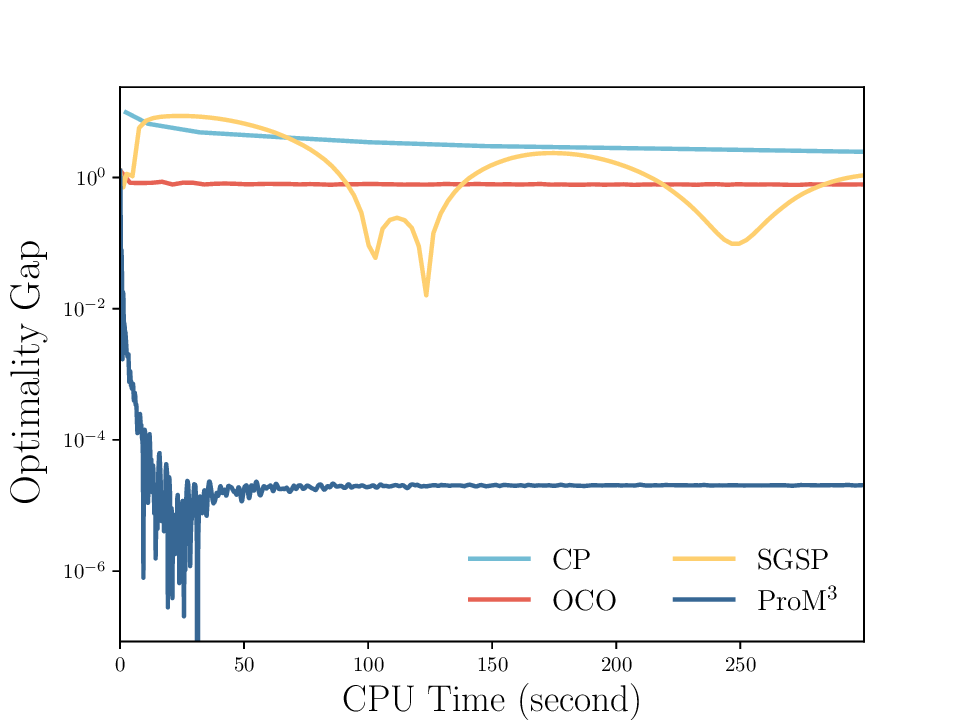}
\includegraphics[scale=0.5]{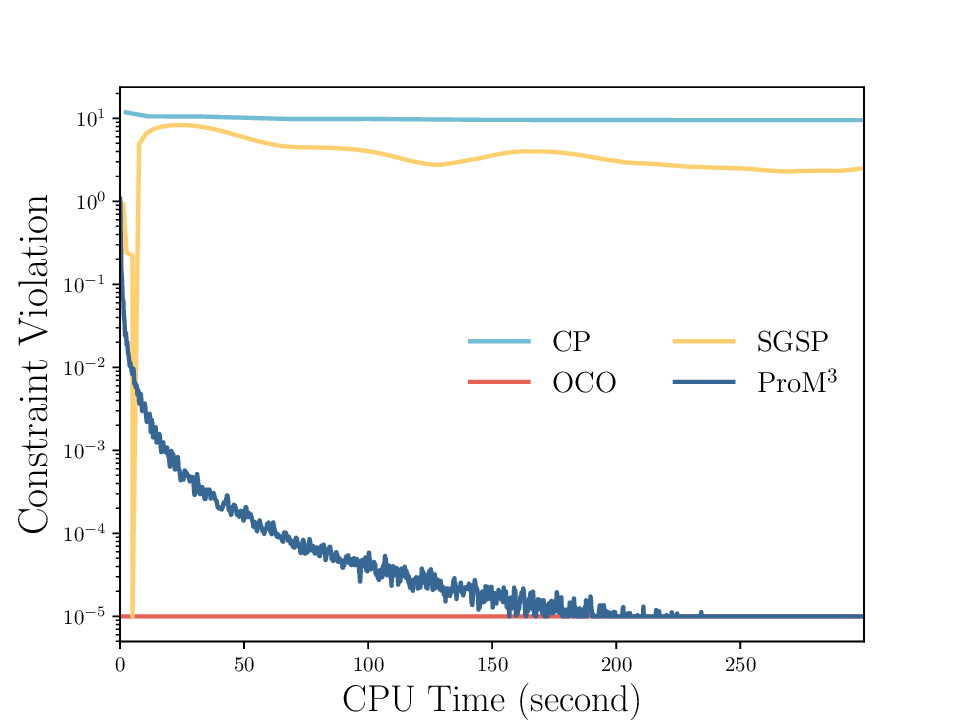}
\end{figure}

\begin{figure}[tb]
\centering
\caption{$(M, N , J_m) = (5, 1000, 200)$. \label{fig-ex2-N=1000_q=200}}
\includegraphics[scale=0.5]{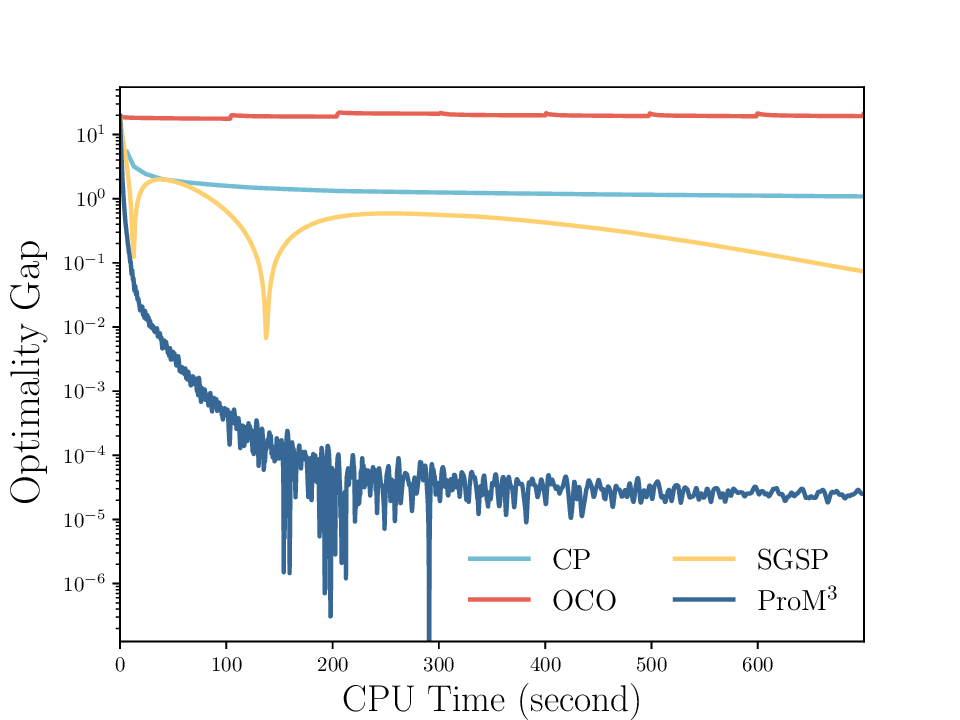}
\includegraphics[scale=0.5]{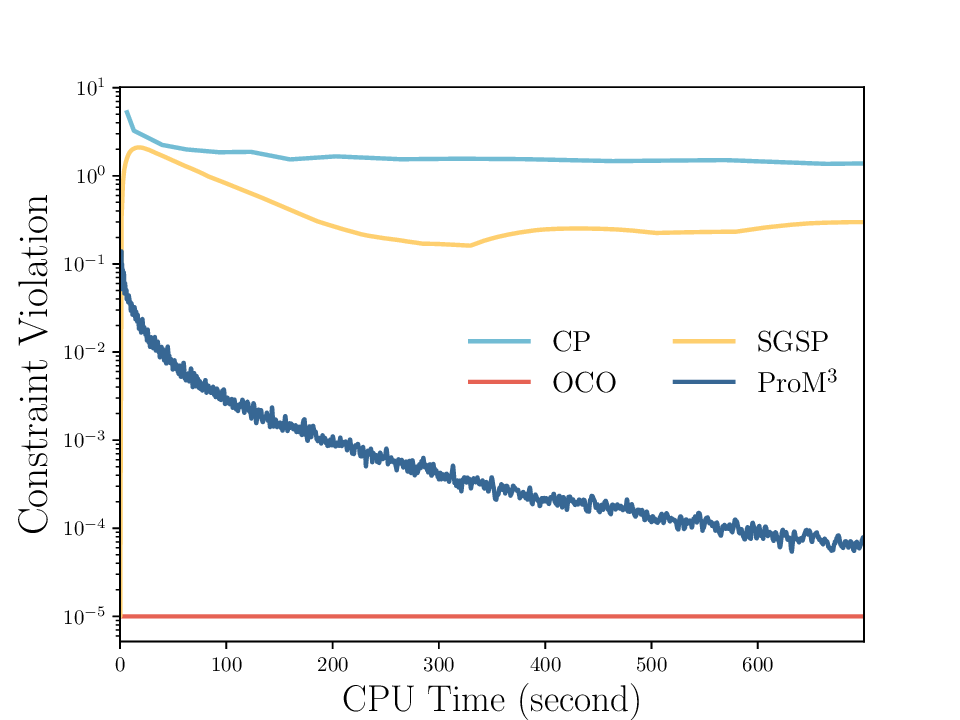}
\end{figure}

\begin{figure}[tb]
\centering
\caption{$(M, N , J_m) = (3, 2000, 200)$. \label{fig-ex2-N=2000_q=200}}
\includegraphics[scale=0.5]{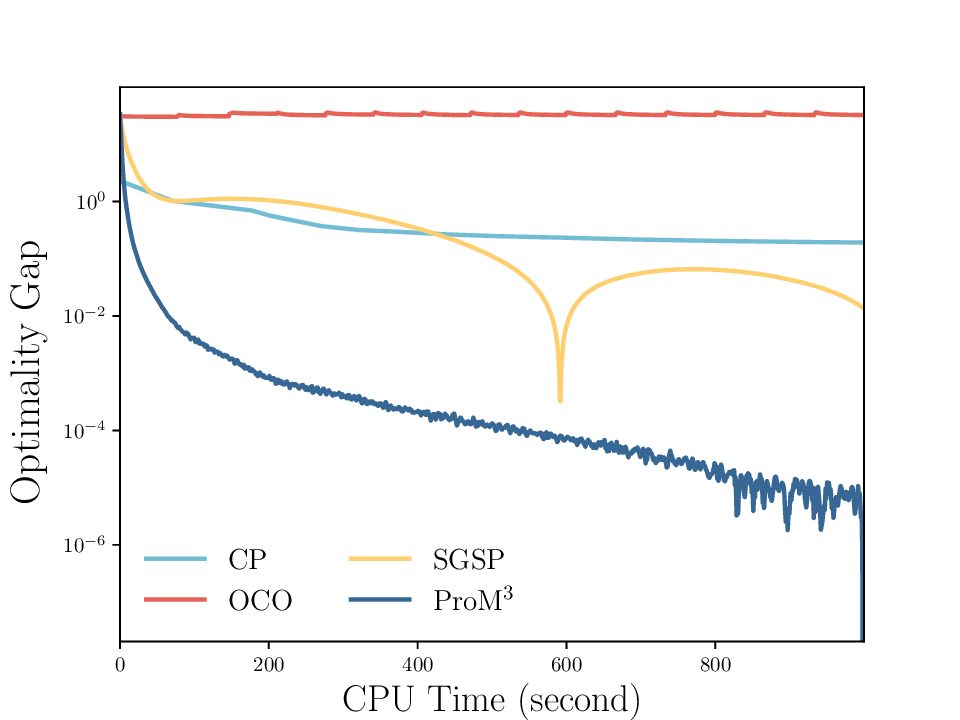}
\includegraphics[scale=0.5]{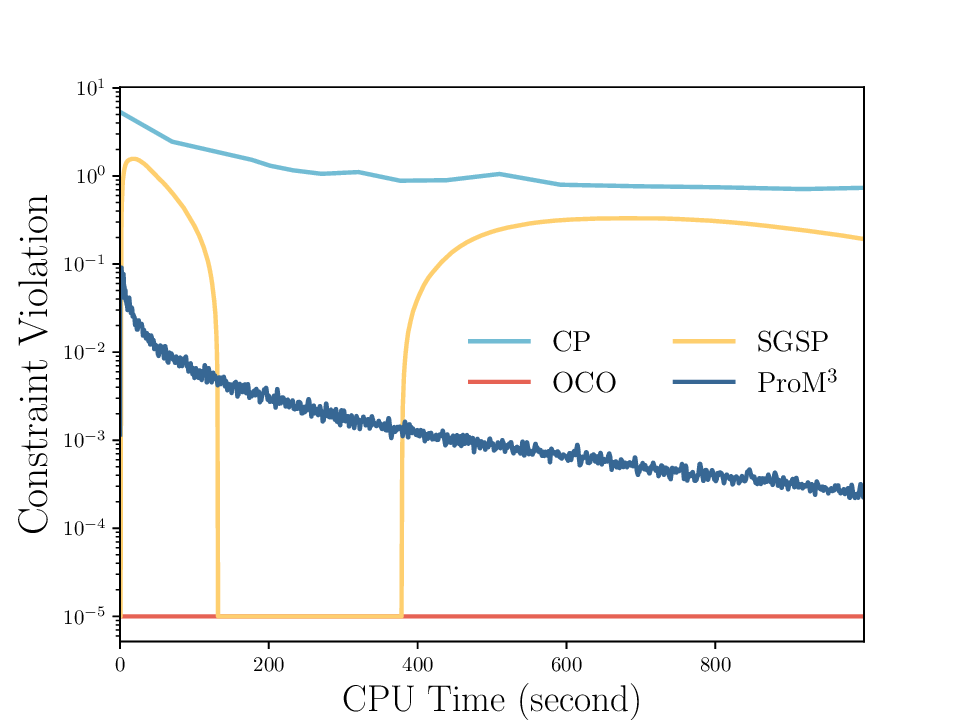}
\end{figure}

The data are generated as follows. We set $\bm{l}=0.001$ and $\bm{u}=1$, and the vector $\bm{c}$ is also generated randomly with i.i.d. standard Gaussian entries. For each $m\in [M]$, denoting $\bm{B}_m = [\bm{b}_{m,2}^{\top};  \bm{b}_{m,3}^{\top} ;\cdots; \bm{b}_{m,J_m}^\top]$, we first generate $\bm{B}_m$ with i.i.d. standard Gaussian entries and then normalize it by $\bm{B}_m \leftarrow \bm{B}_m/\|\bm{B}_m\|_2$. The matrix $\bm{A}_m$ is generated in the same manner as $\bm{B}_m$. We also set $d_m =  \max_{\bm{z}_m \in \mathcal{Z}_m}(\bm{u}/\|\bm{u}\|_2)^{\top} \bm{A}_m \bm{z}_m + \log\left(\bm{z}_{m,1} + \sum_{j=2}^{J_m} \bm{z}_{m,j} \exp{(\bm{b}_{m,j}^{\top}  (\bm{u}/\|\bm{u}\|_2)})\right)$, where $\bm{u}$ is a random vector with entries being i.i.d. uniform random variables on $[0,1]$. 

We test the algorithms on three problem dimensions: $(M, N , J_m) = (5, 200, 1000)$, $(M, N , J_m) = (5, 1000, 200)$ and $(M, N , J_m) = (3, 2000, 200)$, and the corresponding results are plotted in Figures
~\ref{fig-ex2-N=200_q=1000}-\ref{fig-ex2-N=2000_q=200}.
In this experiment, we take the solution returned by our algorithm \prom as the true optimal solution since it achieves a much higher accuracy than the competing methods, as indicated in Section~\ref{sec:robust_QCQP}.

In all these three cases, our algorithm \prom converges much faster in terms of optimality gap than all other algorithms. Also, \prom can achieve the accuracy level $10^{-5}$ to $10^{-6}$, whereas all other methods get stuck at the level of $10^{-1}$. In terms of constraint violation, \prom converges to the feasible region stably and efficiently, whereas the cutting-plane method and SGSP struggle at the level of $10^0$- $10^{-1}$. The sequence of iterates generated by OCO seems to be feasible over the course of execution. This experiment shows that our algorithm is suitable also for highly non-linear RO problems.             

\subsection{Distributionally Robust Newsvendor Problem}

Finally, we investigate the numerical performance of our extend \prom for tackling RO problems with projection-unfriendly uncertainty sets. To this end, we consider a multi-product newsvendor problem that aims at minimizing the total ordering cost subject to a distributionally robust profit-risk constraint.
Let there be $M$ products.
For each product $m\in [M]$, the random demand $d_m$ has $N$ possible outcomes, $d_m^1, \dots, d_m^N$, and the corresponding probabilities are $z_m^1,\dots, z_m^N$. We also denote by $c_m$, $v_m$, $s_m$ and $t_m$ the unit purchase cost, the unit selling price, the unit salvage value and the unit storage cost, respectively. If we purchase $x_m$ units of product $m$, under the demand outcome $d_m^n$, the profit is
\begin{equation}
    \label{eq:profit}
    r( x_m, d_m^n ) = v_m \min\{ d_m^n, x_m \} + s_m ( x_m - d_m^n )_+ - t_m ( d_m^n - x_m )_+ - c_m x_m.
\end{equation} 
The conditional value-at-risk (CVaR) of the loss $- r(x_m, d_m^n)$ at the quantile level $\kappa\in [0,1)$ is given by
\[
\inf_{\tau \in \mathbb{R}} \; \left\{ \frac{\mathbb{E}_{d_m \sim \bm{z}_m}[[ \tau - r(x_m, d_m)]_+] }{1-\kappa } - \tau \right\} ,
\]
which represents the average of the loss $- r(x_m, d_m)$ over its $(1-\kappa)$-tail region. Following the literature on distributionally robust optimization~(\citealt{Ben_Hertog_Waegenaere_Melenberg_Rennen_2013, Esfahani_Kuhn_2017, yue2022linear, gao2023distributionally}), we assume that the probability distribution $\bm{z}_m$ of the random demand is not precisely known but lies in an ambiguity set $\mathcal{Z}_m$. It is then natural to consider the following formulation, which minimizes the purchasing cost subject to the distributionally robust loss CVaR constraint:
\begin{equation}\label{opt:DRNP}
\begin{array}{rll}
\min & \bm{c}^{\top}\bm{x} \\
{\rm s.t.} & \displaystyle \max_{\bm{z}_m \in \mathcal{Z}_m} \left\{ \frac{\mathbb{E}_{ d_m\sim \bm{z}_m }[[ \tau_m - r(x_m, d_m)]_+] }{1-\kappa } - \tau_m \right\} \leq \rho_m  \quad\forall m\in [M] \\
& \bm{x} \in [0,1]^M, \; \bm{\tau}\in \mathbb{R}^M,
\end{array}
\end{equation}
where $\rho_m>0$ is a prescribed risk threshold and $r(x_m, d_m)$ is defined in \eqref{eq:profit}. 
The ambiguity sets $\mathcal{Z}_m$ are taken as an intersection of the probability simplex in $\mathbb{R}^N$ and an Euclidean ball centered at the empirical distribution, \ie, each outcome $d_m^n$ takes probability $1/N$. Formulation~\eqref{opt:DRNP} is thus an instance of RO problems with projection unfriendly uncertainty sets of the form~\eqref{def:generalized-Zm}. The other problem parameters are generated randomly.

\begin{figure}[tb]
\centering
\caption{$(M, N ) = (30, 5000)$. \label{fig-ex3-M=30_N=5000}}
 \includegraphics[scale=0.5]{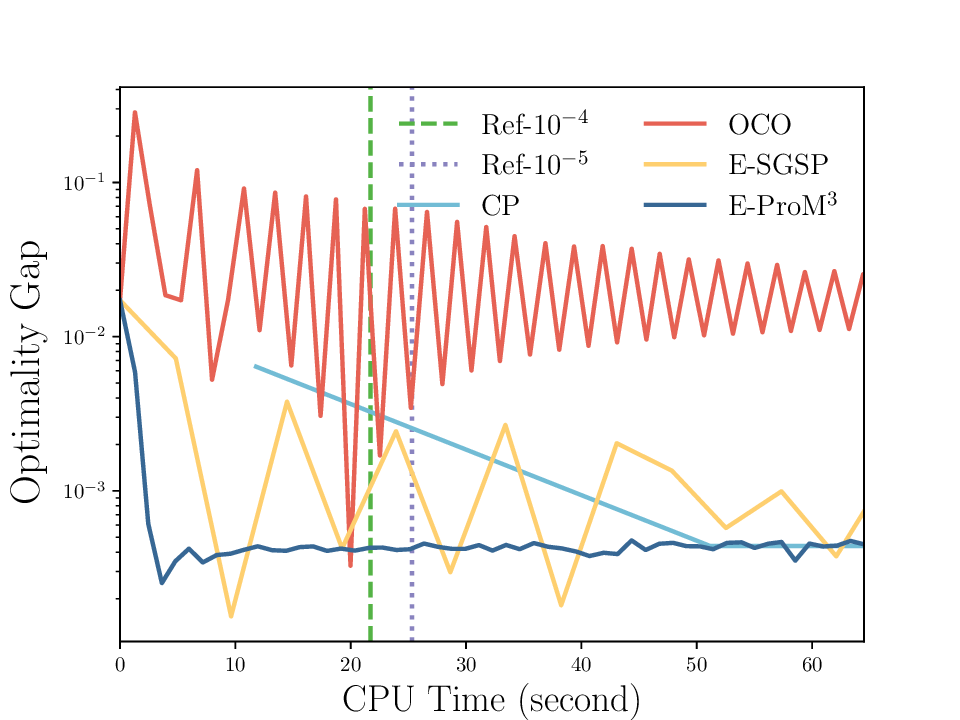}
\includegraphics[scale=0.5]{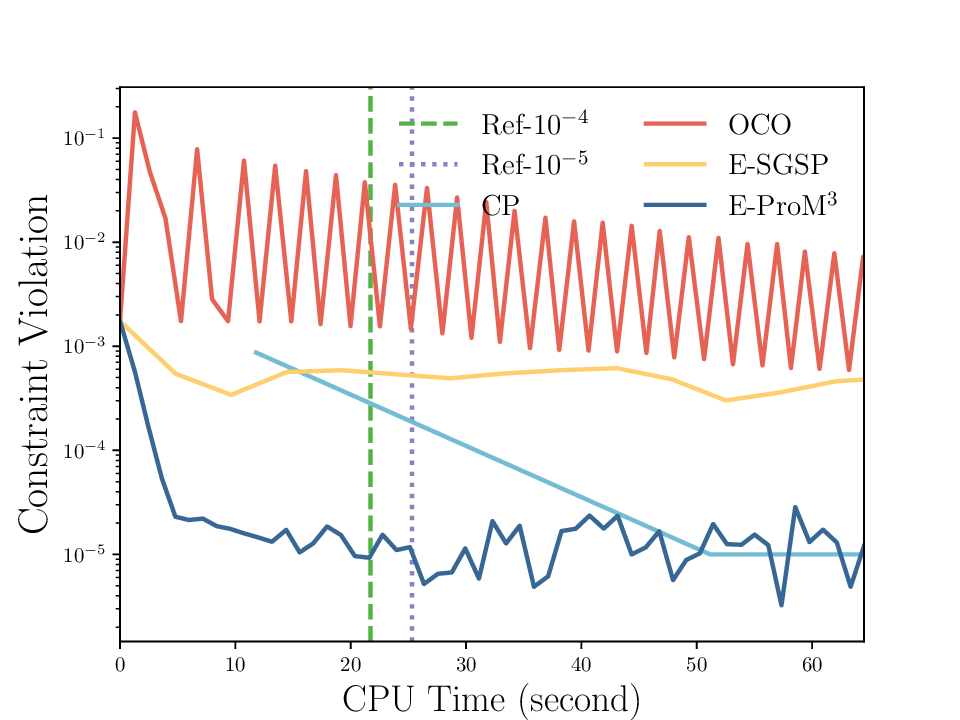}
\end{figure}

Similarly to Sections~\ref{sec:robust_QCQP} and~\ref{sec:log-sum-exp}, we compare our algorithm, the extended \prom, with the cutting-plane method, the reformulation approach and two first-order methods from the papers~\cite{ho2018online} and~\cite{postek2021first}. However, the paper~\cite{postek2021first} also developed an extension of SGSP for RO problems having uncertainty sets of the form \eqref{def:generalized-Zm}. For fairness, we will invoke the extended SGSP in this experiment. The legends ``E-\prom'' and ``E-SGSP'' are adopted to represent the extended \prom and extended SGSP, respectively. For the other competing algorithms, the legends are the same as before. 
The result of a typical instance with dimension $(M,N) = (30, 5000)$ is plotted in Figure~\ref{fig-ex3-M=30_N=5000}. 

The experiment result indicates that our algorithm is more efficient than the competing algorithms, reaching $10^{-3}$ for optimality gap and $10^{-4}$ for constraint violation in a few seconds. The cutting-plane method can also reach the same level of accuracy, for a much longer computational time though. In fact, in this experiment, the cutting-plane method iterated only twice during the 5x seconds of execution. The missing part of the ``CP'' curve at the beginning is due to the fact that the cutting-plane method does not require an initial point $\bm{x}^0$ and an $\bm{x}$-iterate is generated only after the first optimization step is completed.

\section{Conclusion}
Based on a novel max-min-max perspective, this paper devised an iterative algorithm, \prom, for solving RO problems. The algorithm \prom operates directly on the model functions and sets through their gradient and projection oracles, respectively. Such a feature is highly desirable, as it allows for easy exploitation of useful problem structures, and makes the algorithm particularly suitable for contemporary large-scale decision problems. Theoretically, we proved that \prom enjoys strong convergence guarantees. We also extended our algorithm to RO problems with projection-unfriendly uncertainty sets. Numerical results under different challenging regimes (high-dimensional, highly constrained and/or highly non-linear) demonstrated the promising performance of \prom and its extension.

\newpage

\bibliographystyle{ormsv080} 
\bibliography{bibliography} 
%%%%%%%%%%%%%%%%%%%%%%%%%%%%%%%%%%%%%%%%%%%%%%%%%%%%%%%%%%%%%%%%%%%%%%%%%%%%%%%%%%%%
\newpage

\appendix

\section{Auxiliary Results}

\begin{proof}[Proof of \Cref{prop:str_approx_sad_pt_exist}.]
By Corollary~37.6.2 in \cite{rockafellar2015convex}, $\mathcal{F}$ has a saddle point $(\tilde{\bm{u}}, \tilde{\bm{v}}) \in \mathcal{U}\times \mathcal{V}$.
By definition,
\begin{equation}\label{def:saddle-point-F}
\mathcal{F}(\tilde{\bm{u}}, \bm{v}) \leq \mathcal{F}(\tilde{\bm{u}}, \tilde{\bm{v}}) \leq \mathcal{F}(\bm{u}, \tilde{\bm{v}}) \quad\forall (\bm{u},\bm{v})\in \mathcal{U} \times \mathcal{V},
\end{equation}
which implies in particular that $\bm{0}\in \partial_{\bm{u}} \mathcal{F}(\tilde{\bm{u}}, \tilde{\bm{v}})$. By strong convexity of $\mathcal{F}(\cdot, \tilde{\bm{v}})$, we then have 
\[
\mathcal{F}(\bm{u}, \tilde{\bm{v}}) \geq \mathcal{F}(\tilde{\bm{u}}, \tilde{\bm{v}}) +  \frac{\sigma}{2}\|\bm{u} - \tilde{\bm{u}}\|_2^2 ~~~\forall \bm{u} \in \mathcal{U}.
\]
This, together with (\ref{def:saddle-point-F}), implies that
\[
\mathcal{F}(\tilde{\bm{u}}, \bm{v})   \leq \mathcal{F}(\bm{u}, \tilde{\bm{v}}) -\frac{\sigma}{2}\|\bm{u} - \tilde{\bm{u}}\|_2^2 ~~~\forall (\bm{u},\bm{v})\in \mathcal{U} \times \mathcal{V},
\]
concluding the proof.
\end{proof}

\begin{proof}[Proof of \Cref{prop:mmm-equivalence}.]
Each $f_m $ is convex since it is a point-wise supremum of a family of convex functions. Then, by Assumption~\ref{Assump: Robust-Prob-1}\ref{Assump: Robust-Prob-1-iii}-\ref{Assump: Robust-Prob-1-iv} as well as Corollary~28.2.1 and Theorem~28.3 in \cite{rockafellar2015convex}, the max-min problem
\[
\max_{\bm{\lambda} \geq \bm{0}} \min_{\bm{x} \in \mathcal{X}} \; \left\{ f_0(\bm{x}) + \sum_{m \in [M]} \lambda_{m} f_m(\bm{x}) \right\}
\]
admits a saddle point solution, and it is equivalent to the \ref{opt:RO1} problem in the sense that the saddle value equals the optimal value of the \ref{opt:RO1} problem and that for any saddle point $(\bm{\lambda}, \bm{x})$ solving the max-min problem, $\bm{x}$ is an optimal solution to the \ref{opt:RO1} problem.
Unfolding the definitions of $\mathcal{K}$ and $f_m$ completes the proof. 
\end{proof}
\section{Outer Convergence Analysis} 

We first prove that the function $\bm{f}$ is Lipschitz continuous on $\mathcal{X}$.
\begin{proposition}\label{prop:def-lip-f}
Suppose that Assumption~\ref{Assump: Robust-Prob-1}\ref{Assump: Robust-Prob-1-i}-\ref{Assump: Robust-Prob-1-ii} and Assumption~\ref{Assump: Robust-Prob-2} hold.
Then,  
\begin{equation*}
   \| \bm{f}(\bm{x}) - \bm{f} (\bm{y})\|_2 \leq  \mathrm{Lip}_{\bm{f}} \|\bm{x} - \bm{y}\|_2 \quad\forall \bm{x}, \bm{y}\in \mathcal{X} ,
\end{equation*}
where $\mathrm{Lip}_{\bm{f}}=\sqrt{\sum_{m\in [M] }D_m^2}$.
\end{proposition}

\begin{proof}[Proof of \Cref{prop:def-lip-f}.]
By the definition of $f_{m}$, Assumption~\ref{Assump: Robust-Prob-1}\ref{Assump: Robust-Prob-1-i}-\ref{Assump: Robust-Prob-1-ii} and Danskin's theorem (see, \textit{e.g.}, \citealt{bertsekas1997}), we have that for any $m\in [M]$ and $\bm{x} \in \mathcal{X}$, 
\[
\partial f_{m} (\bm{x}) = {\bf conv}\left( \left\{ \bm{\xi}_m \mid \bm{\xi}_m \in  \partial_{\bm{x}} g_{m}(\bm{x}, \bm{z}_{m}^{\star}), ~g_{m}(\bm{x}, \bm{z}_{m}^{\star}) =f_{m}(\bm{x}) \right\} \right),
\]
where ${\bf conv}(\cdot)$ denotes the convex hull.
From Assumption~\ref{Assump: Robust-Prob-2}, any element in $\partial f_{m} (\bm{x}) $ has its norm bounded by $D_m$. By convexity, for any $\bm{x}, \bm{y} \in \mathcal{X}$ and $\hat{\bm{\xi}}_m \in \partial f_{m} (\bm{x})$, 
\[  
f_m(\bm{x}) - f_m(\bm{y}) \le \hat{\bm{\xi}}_m^\top (\bm{x}- \bm{y}) \le D_m \|\bm{x}-\bm{y} \|_2,
\]
which implies 
\[
\|\bm{f}(\bm{x}) - \bm{f}(\bm{y})\|_2 = \sqrt{\sum_{m\in [M] } |f_{m}(\bm{x}) - f_{m} (\bm{y})|^2} \le \sqrt{\sum_{m\in [M] } D_m^2 \|\bm{x}-\bm{y} \|_2^2} = \sqrt{\sum_{m\in [M] } D_m^2 } \|\bm{x}-\bm{y} \|_2,
\]
completing the proof.
\end{proof}

We then prove a technical lemma regarding the saddle-point gap.
\begin{lemma}\label{Prop: K-itertions-property-Algo1}
Suppose that Assumption~\ref{Assump: Robust-Prob-1}\ref{Assump: Robust-Prob-1-i}-\ref{Assump: Robust-Prob-1-ii} and Assumption~\ref{Assump: Robust-Prob-2} hold.
Consider the sequence $\{(\bm{\lambda}^{k}, \bm{x}^{k}, \bm{z}^{k})\}_{k\in [K]}$ generated by Algorithm~\ref{alg:outer} with
$\theta=\frac{1}{K}$, $\nu =\frac{1}{K}$, 
$\alpha \le \frac{1}{\mathrm{Lip}_{\bm{f}}}$ and $\beta \le \frac{1}{2\mathrm{Lip}_{\bm{f}}}$. Then, for any $(\bm{\lambda}, \bm{x}) \in \mathbb{R}^{M}_+ \times \mathcal{X}$, it holds that
\[
\mspace{-20mu}
\sum_{k\in [K]} \left(\mathcal{L}(\bm{\lambda}, \bm{x}^k) - \mathcal{L}(\bm{\lambda}^k, \bm{x})\right) \leq \frac{\|\bm{\lambda} \|_2^2}{2\beta} + \frac{\|\bm{x}-\bm{x}^0\|_2 ^2}{2\alpha} -\frac{\|\bm{\lambda} - \bm{\lambda}^K\|_2^2}{4\beta} +  \frac{3\sqrt{M}}{K} \sum_{k\in [K]} \|\bm{\lambda}^k\|_2 + \left(3 \sqrt{M}\|\bm{\lambda}\|_2 + 1 \right) .
\]
\end{lemma}

\begin{proof}[Proof of \Cref{Prop: K-itertions-property-Algo1}.]
Fix any $(\bm{\lambda}, \bm{x}) \in \mathbb{R}^{M}_+ \times \mathcal{X}$. Since $\mathcal{L}(\bm{\lambda}, \bm{x}) = f_0 (\bm{x}) + \bm{\lambda}^\top \bm{f}(\bm{x})$,
\begin{equation}\label{eq:proof-1}
\begin{split}
&\, \sum_{k\in [K] } \left(\mathcal{L}(\bm{\lambda}, \bm{x}^k) - \mathcal{L}(\bm{\lambda}^k, \bm{x})\right) = \sum_{k = 0}^{K-1} \left(\mathcal{L}(\bm{\lambda}, \bm{x}^{k+1}) - \mathcal{L}(\bm{\lambda}^{k+1}, \bm{x})\right) \\
=&\, \sum_{k = 0}^{K-1} \left( \mathcal{L}(\bm{\lambda}^{k+1},\bm{x}^{k+1}) - \mathcal{L}(\bm{\lambda}^{k+1}, \bm{x}) + \mathcal{L}(\bm{\lambda}, \bm{x}^{k+1}) - \mathcal{L}(\bm{\lambda}^{k+1},\bm{x}^{k+1})   \right) \\
=&\, \sum_{k=0}^{K-1} \left(\mathcal{L}(\bm{\lambda}^{k+1}, \bm{x}^{k+1})- \mathcal{L}(\bm{\lambda}^{k+1}, \bm{x}) + (\bm{\lambda} - \bm{\lambda}^{k+1})^{\top}  \bm{f}( \bm{x}^{k+1})\right).
\end{split}
\end{equation}
Next, for any $\bm{z}\in \mathcal{Z}$, we have
\[
\begin{split}
&\, \mathcal{K}(\bm{\lambda}^{k+1}, \bm{x}^{k+1}, \bm{z})
        - \mathcal{L}(\bm{\lambda}^{k+1}, \bm{x})
\le  \mathcal{K}(\bm{\lambda}^{k+1}, \bm{x}^{k+1}, \bm{z}) - \mathcal{K}(\bm{\lambda}^{k+1}, \bm{x}, \bmt{z}^{k+1}) \\
= &\, f_0( \bm{x}^{k+1} ) + {\bm{\lambda}^{k+1}}^\top \bm{g}( \bm{x}^{k+1}, \bm{z}) + \frac{1}{2\alpha} \|\bm{x}^{k+1} - \bm{x}^k\|_2^2 - f_0( \bm{x} ) - {\bm{\lambda}^{k+1}}^\top \bm{g}( \bm{x}, \tilde{\bm{z}}^{k+1}) - \frac{1}{2\alpha} \|\bm{x}- \bm{x}^k\|_2^2 \\
&\, \quad + \frac{1}{2\alpha} \|\bm{x}- \bm{x}^k\|_2^2 - \frac{1}{2\alpha} \|\bm{x}^{k+1}- \bm{x}^k\|_2^2 \\
\le &\, \frac{1}{2\alpha} \left(\|\bm{x}- \bm{x}^k\|_2^2 - \|\bm{x}- \bm{x}^{k+1}\|_2^2 - \|\bm{x}^k - \bm{x}^{k+1}\|_2^2 \right) + \nu,
\end{split}
\]
where the first inequality (resp., equality) follows from the definition of $\mathcal{L}$ (resp., $\mathcal{K}$) and the second inequality from the fact that $(\bm{x}^{k+1}, \bmt{z}^{k+1})$ is a strong $\nu$-approximate saddle point of the \ref{opt:inner-saddle} problem (see Algorithm~\ref{alg:outer}).
Maximizing the left-hand side w.r.t. $\bm{z}$ over $\mathcal{Z}$ yields
\begin{equation}\label{inequality: step-X-Algo1}
\mathcal{L}(\bm{\lambda}^{k+1}, \bm{x}^{k+1})
- \mathcal{L}(\bm{\lambda}^{k+1}, \bm{x})
\leq \displaystyle \frac{1}{2\alpha} \left(\|\bm{x}- \bm{x}^k\|_2^2 - \|\bm{x}- \bm{x}^{k+1}\|_2^2 - \|\bm{x}^k - \bm{x}^{k+1}\|_2^2 \right) + \nu.
\end{equation}
Also, since
$\bm{\lambda}^{k+1} = \argmax_{\bm{\lambda} \geq \bm{0}} \, \bm{\lambda}^\top \left(2 \bm{g}(\bm{x}^k, \bm{z}^k) - \bm{g}(\bm{x}^{k-1},  \bm{z}^{k-1}) \right) - \frac{1}{2\beta}\|\bm{\lambda} - \bm{\lambda}^k \|_2^2$,
it follows from the optimality of $\bm{\lambda}^{k+1}$ and the strong concavity that 
\[
\frac{1}{2\beta} \left( \|\bm{\lambda}-\bm{\lambda}^{k} \|_2 ^2 -\|\bm{\lambda}-\bm{\lambda}^{k+1} \|_2 ^2 - \|\bm{\lambda}^{k} - \bm{\lambda}^{k+1}\|_2 ^2 \right)  - (\bm{\lambda} - \bm{\lambda}^{k+1})^{\top}   (2 \bm{g}(\bm{x}^{k}, \bm{z}^{k}) -  \bm{g}(\bm{x}^{k-1}, \bm{z}^{k-1})) \geq 0.
\]
Hence,
\begin{equation}\label{inequality: step-lambda-Algo1}
\begin{array}{c@{\;}l}
\displaystyle (\bm{\lambda} - \bm{\lambda}^{k+1})^{\top}  \bm{f}( \bm{x}^{k+1})
\leq& \displaystyle \frac{1}{2\beta} \left( \|\bm{\lambda}-\bm{\lambda}^{k} \|_2 ^2 -\|\bm{\lambda}-\bm{\lambda}^{k+1} \|_2 ^2 - \|\bm{\lambda}^{k} - \bm{\lambda}^{k+1}\|_2 ^2 \right) \\[3mm]
& \displaystyle + (\bm{\lambda} - \bm{\lambda}^{k+1})^{\top} \left( \bm{f}( \bm{x}^{k+1}) - (2 \bm{g}(\bm{x}^{k}, \bm{z}^{k}) -  \bm{g}(\bm{x}^{k-1}, \bm{z}^{k-1})) \right).
\end{array}
\end{equation}
The last term on the right-hand side of~\eqref{inequality: step-lambda-Algo1} satisfies that
\begin{equation}\label{eq:proof-2}
\begin{array}{c@{\;}l}
& \displaystyle (\bm{\lambda} - \bm{\lambda}^{k+1})^{\top} \left(\bm{f}( \bm{x}^{k+1})- (2 \bm{g}(\bm{x}^{k}, \bm{z}^{k}) -  \bm{g}(\bm{x}^{k-1}, \bm{z}^{k-1}) )\right) \\[3mm]
\leq & \displaystyle ( \bm{\lambda} - \bm{\lambda}^{k+1})^{\top} ( \bm{f}( \bm{x}^{k+1}) -  \bm{f}( \bm{x}^{k})) - ( \bm{\lambda} - \bm{\lambda}^{k})^{\top} ( \bm{f}( \bm{x}^{k}) -  \bm{f}( \bm{x}^{k-1})) \\[3mm]
& \displaystyle  +  3\sqrt{M} \theta ( \| \bm{\lambda}\|_2 + \| \bm{\lambda}^{k+1}\|_2) + (\bm{\lambda}^{k+1} - \bm{\lambda}^{k})^{\top} (  \bm{f}( \bm{x}^{k}) -  \bm{f}( \bm{x}^{k-1})) \\[3mm]
\leq & \displaystyle  ( \bm{\lambda} - \bm{\lambda}^{k+1})^{\top} ( \bm{f}( \bm{x}^{k+1}) -  \bm{f}( \bm{x}^{k})) - ( \bm{\lambda} - \bm{\lambda}^{k})^{\top} ( \bm{f}( \bm{x}^{k}) -  \bm{f}( \bm{x}^{k-1})) \\[3mm]
& \displaystyle + 3\sqrt{M} \theta ( \| \bm{\lambda}\|_2 + \| \bm{\lambda}^{k+1}\|_2) + \frac{\mathrm{Lip}_{\bm{f}}}{2}\|\bm{\lambda}^{k+1} -\bm{\lambda}^{k}\|_2 ^2 + \frac{\mathrm{Lip}_{\bm{f}}}{2} \|\bm{x}^{k}- \bm{x}^{k-1}\|_2 ^2,
\end{array}
\end{equation}
where the first inequality follows from the fact that
\[
\| \bm{f}( \bm{x}^k) - \bm{g}(\bm{x}^k, \bm{z}^k) \|_2 =\sqrt{\sum_{m\in [M] } (f_m(\bm{x}^k) - g_m(\bm{x}^k, \bm{z}_m^k) )^2 } \leq \sqrt{M} \theta\quad \forall k\ge 0,
\]
and the second inequality from \Cref{prop:def-lip-f}.
Substituting~\eqref{eq:proof-2} into~\eqref{inequality: step-lambda-Algo1}, we get
\begin{equation*}\label{eq:proof-3}
\begin{split}
(\bm{\lambda} - \bm{\lambda}^{k+1})^{\top}  \bm{f}( \bm{x}^{k+1})
\le &\, \frac{1}{2\beta} \left( \|\bm{\lambda}-\bm{\lambda}^{k} \|_2 ^2 -\|\bm{\lambda}-\bm{\lambda}^{k+1} \|_2 ^2 - \|\bm{\lambda}^{k} - \bm{\lambda}^{k+1}\|_2 ^2 \right) \\
&\,  + ( \bm{\lambda} - \bm{\lambda}^{k+1})^{\top} ( \bm{f}( \bm{x}^{k+1}) -  \bm{f}( \bm{x}^{k})) - ( \bm{\lambda} - \bm{\lambda}^{k})^{\top} ( \bm{f}( \bm{x}^{k}) -  \bm{f}( \bm{x}^{k-1})) \\
&\,  + 3\sqrt{M} \theta ( \| \bm{\lambda}\|_2 + \| \bm{\lambda}^{k+1}\|_2) + \frac{\mathrm{Lip}_{\bm{f}}}{2}\|\bm{\lambda}^{k+1} -\bm{\lambda}^{k}\|_2 ^2 + \frac{\mathrm{Lip}_{\bm{f}}}{2} \|\bm{x}^{k}- \bm{x}^{k-1}\|_2 ^2,
\end{split}
\end{equation*}
Then, substituting the above inequality and inequality~\eqref{inequality: step-X-Algo1} into \eqref{eq:proof-1},  we obtain
\[
\mspace{-20mu}
\begin{array}{l@{\;}l}
& \displaystyle \sum_{k\in [K] } \left(\mathcal{L}(\bm{\lambda}, \bm{x}^{k}) - \mathcal{L}( \bm{\lambda}^{k}, \bm{x})\right) \\
\leq& \displaystyle \frac{\| \bm{x}-\bm{x}^{0}\|_2 ^2}{2\alpha} +\frac{\| \bm{\lambda} \|_2 ^2}{2\beta} - \frac{ \| \bm{x}- \bm{x}^{K}\|_2 ^2 }{2\alpha} - \frac{ \| \bm{\lambda} - \bm{\lambda}^{K}\|_2 ^2}{2\beta} - \frac{\|\bm{x}^{K}-\bm{x}^{K-1}\|^2}{2\alpha}    \\
& \displaystyle + (\bm{\lambda} - \bm{\lambda}^{K})^{\top} ( \bm{f}( \bm{x}^{K}) -  \bm{f}( \bm{x}^{K-1})) + 3\sqrt{M} \theta  \sum _{k\in [K] } \|\bm{\lambda}^{k}\|_2 + 3\sqrt{M} \theta K \|\bm{\lambda}\|_2 + \nu K\\
&\displaystyle  -\left( \frac{1}{2\beta} -\frac{\mathrm{Lip}_{\bm{f}}}{2}\right)\sum_{k\in [K]}\|\bm{\lambda}^{k} -\bm{\lambda}^{k-1}\|_2^2 - \left(\frac{1}{2\alpha} - \frac{\mathrm{Lip}_{\bm{f}}}{2} \right) \sum_{k\in [K-1]}\|\bm{x}^k - \bm{x}^{k-1}\|_2^2  \\
\leq &\displaystyle \frac{\|\bm{x}-\bm{x}^{0}\|_2 ^2}{2\alpha}  + \frac{\|\bm{\lambda}\|_2 ^2}{2\beta} - \frac{\|\bm{x}- \bm{x}^{K}\|_2 ^2 }{2\alpha} - \frac{\|\bm{\lambda} - \bm{\lambda}^{K} \|_2^2}{4\beta}+ 3 \sqrt{M} \theta \sum _{k\in [K] } \|\bm{\lambda}^{k}\|_2+ \left(3\sqrt{M} \theta \|\bm{\lambda}\|_2 + \nu \right) K,   \\
& \displaystyle  - \left(\frac{1}{4\beta} -\frac{\mathrm{Lip}_{\bm{f}}}{2}\right) \|\bm{\lambda} - \bm{\lambda}^{K}\|_2^2 - \left(\frac{1}{2\alpha} - \frac{\mathrm{Lip}_{\bm{f}}}{2}\right) 
\sum_{k\in [K]}\|\bm{x}^{k}- \bm{x}^{k-1}\|^2_2 -\left( \frac{1}{2\beta} -\frac{\mathrm{Lip}_{\bm{f}}}{2}\right)\sum_{k\in [K]}\|\bm{\lambda}^{k} -\bm{\lambda}^{k-1}\|_2^2 \\
\leq & \displaystyle \frac{\|\bm{x}-\bm{x}^{0}\|_2 ^2}{2\alpha}  + \frac{\|\bm{\lambda}\|_2 ^2}{2\beta}  - \frac{\|\bm{\lambda} - \bm{\lambda}^{K} \|_2^2}{4\beta}+ \frac{3 \sqrt{M}}{K} \sum _{k\in [K] } \|\bm{\lambda}^{k}\|_2+ \left(3\sqrt{M}  \|\bm{\lambda}\|_2 + 1 \right) ,
\end{array}
\]
where the first inequality follows from telescoping and the initial conditions $\bm{\lambda}^0 = \bm{0}$ and $\bm{x}^{-1}=\bm{x}^{0}$, and the second from the inequality
\[
(\bm{\lambda} - \bm{\lambda}^{K})^{\top} ( \bm{f}( \bm{x}^{K}) -  \bm{f}( \bm{x}^{K-1})) \leq \frac{\mathrm{Lip}_{\bm{f}}}{2}( \|\bm{\lambda} - \bm{\lambda}^{K}\|_2 ^2 + \|\bm{x}^{K}- \bm{x}^{K-1}\|_2 ^2),
\]
and the third from the choice of $\theta$, $\nu$, $\alpha$ and $\beta$. This completes the proof.
\end{proof}

The following lemma asserts that the iterator $\bm{\lambda}^k$ is upper bounded uniformly in $k$.

\begin{lemma}\label{Prop: lambda-Upper-bound-Algo1}
Suppose that Assumptions~\ref{Assump: Robust-Prob-1} and \ref{Assump: Robust-Prob-2} hold. Then, $\mathcal{L}$ admits a saddle point. 
Moreover, the sequence $\{\bm{\lambda}^k\}_{k\in [K]}$ generated by the Algorithm~\ref{alg:outer} with $\theta=\frac{1}{K}$, $\nu =\frac{1}{K}$, 
$\alpha \le \frac{1}{\mathrm{Lip}_{\bm{f}}}$ and $\beta \le \frac{1}{2\mathrm{Lip}_{\bm{f}}}$ satisfies that for any $k\in [K]$, 
\begin{align*} 
% \label{inequality: lambda-Upper-bound-Algo1}
\|\bm{\lambda}^{k}\|^2 _2\le 12\|\bm{\lambda}^{\star}\|_2^2 +\frac{8\beta}{\alpha} \|\bm{x}^{\star} - \bm{x}^{0}\|^2_2 +48 \beta \sqrt{M}  \|\bm{\lambda}^{\star}\|_2   + \frac{576\beta M}{\mathrm{Lip}_{\bm{f}}}   + 16 \beta ,
\end{align*}
where $(\bm{\lambda}^\star, \bm{x}^\star) \in \mathbb{R}^M_+\times \mathcal{X}$ is an arbitrary saddle point of $\mathcal{L}$.
\end{lemma}

\begin{proof}[Proof of \Cref{Prop: lambda-Upper-bound-Algo1}.]
The existence of a saddle point of $\mathcal{L}$ follows from the proof of \Cref{prop:mmm-equivalence}.
By definition, for any $(\bm{\lambda}, \bm{x}) \in \mathbb{R}^{M}_{+} \times \mathcal{X}$, we have
\begin{equation}\label{def: L-Saddle-Point}
\mathcal{L}(\bm{\lambda}, \bm{x}^\star) \leq \mathcal{L}(\bm{\lambda}^\star, \bm{x}^\star) \leq \mathcal{L}(\bm{\lambda}^\star, \bm{x}).
\end{equation}
Taking $(\bm{\lambda}, \bm{x}) = (\bm{\lambda}^k, \bm{x}^k)$, we then have $\mathcal{L}(\bm{\lambda}^{\star}, \bm{x}^k) \geq \mathcal{L}(\bm{\lambda}^k, \bm{x}^{\star}) ~\forall k \in [K]$, which implies that
\begin{equation}\label{inequality: lambda-Upper-bound-1}
\displaystyle \frac{\|\bm{\lambda}^{\star} - \bm{\lambda}^{k'}\|_2^2}{4\beta} \leq\displaystyle \sum _{k\in [k'] } \left( \mathcal{L}(\bm{\lambda}^{\star}, \bm{x}^{k}) - \mathcal{L}(\bm{\lambda}^{k}, \bm{x}^{\star})\right) + \frac{\|\bm{\lambda}^{\star} - \bm{\lambda}^{k'}\|_2^2}{4\beta} ~\forall k' \in [K]. 
\end{equation}
Using \Cref{Prop: K-itertions-property-Algo1} with $(\bm{\lambda}, \bm{x})=(\bm{\lambda}^{\star}, \bm{x}^{\star})$ and \eqref{inequality: lambda-Upper-bound-1}, we obtain 
\[
\begin{array}{l@{\;}c@{\;}l}
\displaystyle \frac{\|\bm{\lambda}^{\star} - \bm{\lambda}^{k'}\|_2^2}{4\beta} 
&\leq& \displaystyle \frac{\|\bm{\lambda}^{\star}\|_2 ^2}{2\beta} + \frac{\|\bm{x}^{\star}-\bm{x}^{0}\|_2 ^2}{2\alpha} + 3 \sqrt{M} \theta  \sum _{k\in [k'] } \|\bm{\lambda}^{k}\|_2+ 3 \sqrt{M}  \|\bm{\lambda}^{\star}\|_2  + 1 \\
&\leq& \displaystyle \frac{\|\bm{\lambda}^{\star}\|_2 ^2}{2\beta} + \frac{\|\bm{x}^{\star}-\bm{x}^{0}\|_2 ^2}{2\alpha} +  \sum_{k\in [k']} \left( \frac{\|\bm{\lambda}^{k}\|_2^2 \mathrm{Lip}_{\bm{f}}}{16K } + \frac{36 K M \theta^2}{\mathrm{Lip}_{\bm{f}}} \right)  + 3 \sqrt{M}  \|\bm{\lambda}^{\star}\|_2  + 1 \\
&\leq & \displaystyle \frac{\|\bm{\lambda}^{\star}\|_2^2}{2\beta} + \frac{\|\bm{x}^{\star}-\bm{x}^{0}\|_2^2}{2\alpha} + \frac{ \mathrm{Lip}_{\bm{f}}}{16 K} \sum _{k\in [k'] }\|\bm{\lambda}^{k}\|_2 ^{2} + \frac{36 M}{\mathrm{Lip}_{\bm{f}}} +3  \sqrt{M} \|\bm{\lambda}^{\star}\|_2  + 1,
\end{array}
\]
where the last inequality follows from the fact that $\theta \le \frac{1}{K}$.
Since $\beta \le \frac{1}{2\mathrm{Lip}_{\bm{f}} }$, 
\begin{equation*}
    \|\bm{\lambda}^{\star} - \bm{\lambda}^{k'}\|_2^2 \le 2\|\bm{\lambda}^{\star}\|_2^2 + \frac{2\beta\|\bm{x}^{\star}-\bm{x}^{0}\|_2^2}{\alpha} + \frac{ 1}{8 K} \sum _{k\in [k'] }\|\bm{\lambda}^{k}\|_2 ^{2} + \frac{144 \beta M}{\mathrm{Lip}_{\bm{f}}} + 12\beta \sqrt{M} \|\bm{\lambda}^{\star}\|_2  + 4\beta.
\end{equation*}
For simplicity, we denote
\begin{equation}
    \label{def:Gamma}
    \Gamma = 12\|\bm{\lambda}^{\star}\|_2^2 +\frac{8\beta}{\alpha} \|\bm{x}^{\star} - \bm{x}^{0}\|^2_2 +48 \beta \sqrt{M}  \|\bm{\lambda}^{\star}\|_2   + \frac{576\beta M}{\mathrm{Lip}_{\bm{f}}}   + 16 \beta.
\end{equation}
Then, the above inequality implies that for any $k'\in [K]$,
\[
\|\bm{\lambda}^{\star} - \bm{\lambda}^{k'}\|_2^2 \leq \displaystyle \frac{1}{4} \Gamma - \|\bm{\lambda}^{\star}\|_2^2 +\frac{1}{8K} \sum_{k\in [K] }\|\bm{\lambda}^{k}\|_2^2.
\]
This, together with the fact that $\|\bm{\lambda}^{k'}\|_2^2 \leq 2\|\bm{\lambda}^{\star} - \bm{\lambda}^{k'}\|^2_2 + 2 \|\bm{\lambda}^{\star}\|^2_2$ for all $k' \in [K]$, yields
\begin{equation}\label{inequality: lambda-Upper-bound-2}
\|\bm{\lambda}^{k'}\|_2^2 \leq \displaystyle  \frac{1}{2}	\Gamma +\frac{ 1}{4K} \sum_{k\in [K] }\|\bm{\lambda}^{k}\|_2^2 ~~~\forall k' \in [K].
\end{equation}
% By induction, we now establish the upper bound~\eqref{inequality: lambda-Upper-bound-Algo1}. 
Next, we prove by induction the desired result that $\|\bm{\lambda}^{k'}\|_2 \le \Gamma$ for any $k'\in [K]$. 
For $k' =1$, by \eqref{inequality: lambda-Upper-bound-2}, we have 
\[
\|\bm{\lambda}^{1}\|_2^2 \leq \frac{1}{2}	\Gamma + \frac{1}{4} \|\bm{\lambda}^{1}\|_2^2\leq \frac{1}{2}	\Gamma + \frac{1}{2} \|\bm{\lambda}^{1}\|_2^2,
\]
which implies
\[
\|\bm{\lambda}^{1}\|_2^2 \leq \Gamma.
\]
Next, assume $\|\bm{\lambda}^{k'}\|_2 \le \Gamma$ holds for some $k' \in [1, K)$. By~\eqref{inequality: lambda-Upper-bound-2} we then have 
\[
\begin{array}{l@{\;}l}
\displaystyle \|\bm{\lambda}^{k'+1}\|_2^2 &\leq \displaystyle \frac{1}{2} \Gamma + \frac{1}{4K}\sum_{k=1}^{k'+1} \|\bm{\lambda}^{k}\|_2^2 =\frac{1}{2}	 \Gamma   + \frac{1}{4}	\|\bm{\lambda}^{k'+1}\|_2^2 + \frac{1}{4K}\sum_{k=1}^{k'} \|\bm{\lambda}^{k}\|_2^2 \\
&\leq \displaystyle \frac{1}{2}	\Gamma  + \frac{1}{4}	\|\bm{\lambda}^{k'+1}\|_2^2 + \frac{1}{4K}\sum_{k=1}^{k'} \Gamma \leq \frac{3}{4} \Gamma + \frac{1}{4}	\|\bm{\lambda}^{k'+1}\|_2^2,
\end{array}
\]
which concludes our proof.
\end{proof}

We are now ready to prove \Cref{thm: iteration-complexity-Algo1}.

\begin{proof}[Proof of \Cref{thm: iteration-complexity-Algo1}.]
To bound $f_0(\bar{\bm{x}}^K) - f_0(\bm{x}^{\star})$, we let $\bar{\bm{\lambda}}^K=\frac{1}{K} \sum_{k\in [K] }\bm{\lambda}^{k}$. By the convexity of $\mathcal{L}(\bm{\lambda}, \cdot)$ and $- \mathcal{L}(\cdot, \bm{x})$, and Lemmas~\ref{Prop: K-itertions-property-Algo1} and~\ref{Prop: lambda-Upper-bound-Algo1}, we have that for any $(\bm{\lambda},\bm{x}) \in \mathbb{R}^M_+ \times \mathcal{X}$,
\begin{align}\label{inequality: thm-iteration-comp-2}
\mathcal{L}(\bm{\lambda}, \bar{\bm{x}}^K) - \mathcal{L}(\bar{\bm{\lambda}}^K, \bm{x})  \leq \frac{1}{K} \left( \frac{\|\bm{\lambda}\|_2^2}{2\beta} + \frac{\|\bm{x}-\bm{x}^0\|_2^2}{2\alpha} + 3 \sqrt{M \Gamma }  + 3 \sqrt{M}\|\bm{\lambda}\|_2   +  1 \right) ,
\end{align}
where $\Gamma$ is the constant defined in \eqref{def:Gamma}. Then,
\begin{equation*}
\begin{split}
&\,f_0(\bar{\bm{x}}^K) - f_0(\bm{x}^{\star}) =  \mathcal{L}(\bm{0}, \bar{\bm{x}}^K) - \mathcal{L}( \bm{\lambda}^{\star}, \bm{x}^{\star}) \leq \mathcal{L}( \bm{0}, \bar{\bm{x}}^K) - \mathcal{L}( \bar{\bm{\lambda}}, \bm{x}^{\star})\\
\le &\, \frac{1}{K} \left(  \frac{\|\bm{x}^\star -\bm{x}^0\|_2^2}{2\alpha} + 3 \sqrt{M \Gamma }   +  1 \right),
\end{split}
\end{equation*}
where the equality follows from the definition of $\mathcal{L}$, the first inequality from \eqref{def: L-Saddle-Point}, and the second inequality from \eqref{inequality: thm-iteration-comp-2}.

Next, we bound $\max_{m\in [M]}  [f_{m}(\bar{\bm{x}}^K)]_{+} $. By the definition of $\mathcal{L}$ and \eqref{def: L-Saddle-Point}, we have
\begin{equation}
    \label{inequality: non-negative-f-0}
    f_0(\bar{\bm{x}}^K) - f_0(\bm{x}^{\star}) \geq - \sum_{m\in [M] }\lambda_{m}^{\star} f_{m}(\bar{\bm{x}}^K) \geq  - \sum_{m\in [M] }\lambda_{m}^{\star} [f_{m}(\bar{\bm{x}}^K)]_{+}.
\end{equation}
Let $\hat{\bm{\lambda}}$ be the vector defined by $\hat{\lambda}_m=1+ \lambda_m^{\star}$ if $f_m(\bar{\bm{x}}^{K})>0$, and $\hat{\lambda}_m=0$ otherwise, for any $m\in [M]$. 
By the definition of $\mathcal{L}$, \eqref{def: L-Saddle-Point} and \eqref{inequality: thm-iteration-comp-2}, we have 
\begin{align*}
    &\, f_0(\bar{\bm{x}}^K ) - f_0(\bm{x}^{\star})  + \sum_{m\in [M]} \hat{\lambda}_m f_m(\bar{\bm{x}}^K)  = \mathcal{L}(\hat{\bm{\lambda}},\bar{\bm{x}}^K) - \mathcal{L}(\bm{\lambda}^{\star},\bm{x}^{\star}) \\
    \le &\, \mathcal{L}(\hat{\bm{\lambda}}, \bar{\bm{x}}^K) - \mathcal{L}(\bar{\bm{\lambda}}^K, \bm{x}^\star) \le  \frac{1}{K} \left( \frac{\|\hat{\bm{\lambda}}\|_2^2}{2\beta} + \frac{\|\bm{x}^\star - \bm{x}^0\|_2^2}{2\alpha} + 3 \sqrt{M \Gamma }  + 3 \sqrt{M}\|\hat{\bm{\lambda}}\|_2   +  1 \right) ,
\end{align*}
which, together with \eqref{inequality: non-negative-f-0}, yields
\begin{align*}
    &\, \max_{m\in [M]}  [f_{m}(\bar{\bm{x}}^K )]_{+}  \le   \sum_{m\in [M] }  [f_{m}(\bar{\bm{x}}^K)]_{+} + \sum_{m\in [M] }\lambda_{m}^{\star} [f_{m}(\bar{\bm{x}}^K)]_{+} + f_0(\bar{\bm{x}}^K ) - f_0(\bm{x}^{\star})\\
    = &\, f_0(\bar{\bm{x}}^K ) - f_0(\bm{x}^{\star})  + \sum_{m\in [M]} \hat{\lambda}_m f_m(\bar{\bm{x}}^K) \le  \frac{1}{K} \left( \frac{\|\hat{\bm{\lambda}}\|_2^2}{2\beta} + \frac{\|\bm{x}^\star - \bm{x}^0\|_2^2}{2\alpha} + 3 \sqrt{M \Gamma }  + 3 \sqrt{M}\|\hat{\bm{\lambda}}\|_2   +  1 \right).
\end{align*}
This completes the proof.
\end{proof}

\section{Inner Convergence Analysis} 
\label{sec:inner-conv-analysis}

To analyze the \ref{opt:inner-saddle} problem, we consider general saddle-point problems of the form
\begin{equation}\label{opt:inner-saddle-general}
\min_{\bm{u}\in \mathcal{U}} \; \max_{\bm{v}\in \mathcal{V}} \; \hat{\mathcal{F}}(\bm{u}, \bm{v}) + \frac{\sigma}{2}\| \bm{u} - \hat{\bm{u}} \|_2^2,
\end{equation}
where $\hat{\mathcal{F}}(\bm{u}, \bm{v})$ is convex in $\bm{u}$ and concave in $\bm{v}$ (possibly non-linear, non-smooth), $\mathcal{U}$ and $\mathcal{V}$ are non-empty compact convex sets, and $\hat{\bm{u}}\in\mathcal{U}$. We solve it by the following algorithm.
\begin{algorithm}[htbp]
\caption{Modified Projected Subgradient Ascent Descent. \label{alg:mpsad}}
\SetKwInOut{Input}{Input}
\SetKwInOut{Output}{Output}
\Input{$T\ge 1$,  $\delta >0$,  $\gamma>0$, $\bm{u}_{0} \in \mathcal{U}$, $\bm{v}_0\in \mathcal{V}$.}
\For{$t = 0, \dots, T-1$}
{
($\bm{v}$-update) Compute $\bm{\zeta}_t \in \partial_{\bm{v}} (- \hat{\mathcal{F}})( \bm{u}_t, \bm{v}_t )$. Set
\[ 
\bm{v}_{t+1} = \Proj_{\mathcal{V}} \left( \bm{v}_t - \delta \, (2\bm{\zeta}_t - \bm{\zeta}_{t-1} ) \right). 
\]

($\bm{u}$-update) Compute $\bm{\xi}_t \in \partial_{\bm{u}} \hat{\mathcal{F}} (\bm{u}_t, \bm{v}_{t+1})$. Set 
\[
\bm{u}_{t+1}  =  \Proj_{\mathcal{U}} \left( \frac{1}{1 + \gamma \sigma }\left(\gamma \sigma  \hat{\bm{u}} +  \bm{u}_t - \gamma \bm{\xi}_t\right) \right).
\]
}
\Output{$\bar{\bm{u}}_T =\frac{1}{T} \sum_{t\in [T]}\bm{u}_{t} $ and $\bar{\bm{v} }_T = \frac{1}{T}\sum_{t\in [T]}\bm{v}_{t}$.}
\end{algorithm}

\begin{proposition}\label{thm: iteration-complexity-Algo2}
Suppose that $\mathcal{U}$ and $\mathcal{V}$ are nonempty compact convex sets, that $\hat{\mathcal{F}}(\cdot, \bm{v})$ is convex on $\mathcal{U}$ for any $\bm{v}\in \mathcal{V}$ and $\hat{\mathcal{F}}(\bm{u}, \cdot)$ is  concave on $\mathcal{V}$ for any $\bm{u}\in \mathcal{U}$, and that there exist $C_1, C_2>0$ such that $\|\bm{\xi}\|_2 \leq C_1$ and $\|\bm{\zeta}\|_2\leq C_2$ for any $(\bm{u}, \bm{v}) \in \mathcal{U} \times \mathcal{V}$, $\bm{\xi} \in \partial_{\bm{u} } \hat{\mathcal{F}}(\bm{u}, \bm{v})$ and $\bm{\zeta}\in \partial_{\bm{v}}(-\hat{\mathcal{F}})(\bm{u}, \bm{v})$.
Then, the output $(\bar{\bm{u}}_T, \bar{\bm{v}}_T)$ of Algorithm~\ref{alg:mpsad} with $\gamma , \delta \le \tfrac{1}{\sqrt{T}}$ is a strong $\mathcal{O}(T^{-1/2})$-approximate saddle point of problem~\eqref{opt:inner-saddle-general}. 
\end{proposition}

\begin{proof}[Proof of \Cref{thm: iteration-complexity-Algo2}.]
We first note that
\begin{equation}\label{equality: thm-iter-comp-algo2-1}
    \begin{split}
        \hat{\mathcal{F}}(\bm{u}_{t+1}, \bm{v}) - \hat{\mathcal{F}}(\bm{u}, \bm{v}_{t+1}) &= \hat{\mathcal{F}}(\bm{u} _{t+1}, \bm{v})  -\hat{\mathcal{F}}(\bm{u} _{t+1}, \bm{v}_{t+1}) + \hat{\mathcal{F}}(\bm{u}_{t+1}, \bm{v}_{t+1}) - \hat{\mathcal{F}}(\bm{u}_t, \bm{v}_{t+1}) \\
        & + \hat{\mathcal{F}}(\bm{u}_{t}, \bm{v}_{t+1}) - \hat{\mathcal{F}}( \bm{u}, \bm{v}_{t+1}).
    \end{split}
\end{equation}
Since $\hat{\mathcal{F}}(\cdot, \bm{v})$ is convex on $\mathcal{U}$ for any $\bm{v}\in \mathcal{V}$ and $\hat{\mathcal{F}}(\bm{u}, \cdot)$ is  concave on $\mathcal{V}$ for any $\bm{u}\in \mathcal{U}$, for any $\bm{\xi}'_{t}\in \partial_{\bm{u}} \hat{\mathcal{F}}( \bm{u} _{t+1}, \bm{v}_{t+1}) $ and $\bm{\zeta}_{t+1} \in \partial_{\bm{v}} (-\hat{\mathcal{F}})(\bm{u}_{t+1}, \bm{v}_{t+1})$, we get the three inequalities
\begin{equation}\label{eq:proof-5}
    \begin{split}
        \hat{\mathcal{F}}(\bm{u}_{t+1}, \bm{v}) - \hat{\mathcal{F}}(\bm{u}_{t+1}, \bm{v}_{t+1}) &\leq -(\bm{v} - \bm{v}_{t+1})^{\top} \bm{\zeta}_{t+1}\\
\hat{\mathcal{F}}(\bm{u}_{t+1}, \bm{v}_{t+1}) - \hat{\mathcal{F}}(\bm{u}_t, \bm{v}_{t+1}) &\leq - (\bm{u}_{t} - \bm{u}_{t+1})^{\top}\bm{\xi}'_t\\
\hat{\mathcal{F}}(\bm{u}_{t}, \bm{v}_{t+1}) - \hat{\mathcal{F}}(\bm{u}, \bm{v}_{t+1}) &\leq - (\bm{u} -\bm{u} _t )^{\top} \bm{\xi}_{t}
    \end{split}
\end{equation}
Noting that 
\begin{align*}
    \bm{u}_{t+1} & = \argmin_{\bm{u}\in \mathcal{U}}\; \bm{\xi}_t^\top (\bm{u} - \bm{u}_t) + \tfrac{\sigma}{2}\| \bm{u} - \hat{\bm{u}} \|_2^2 + \tfrac{1}{2\gamma}\| \bm{u} - \bm{u}_t \|_2^2 \\
    \bm{v}_{t+1} & = \argmin_{\bm{v}\in \mathcal{V}}\; (2\bm{\zeta}_t - \bm{\zeta}_{t-1} )^\top (\bm{v} - \bm{v}_t)  + \tfrac{1}{2\delta}\| \bm{v} - \bm{v}_t \|_2^2
\end{align*}
and using their optimality conditions, we get the two inequalities
\begin{equation}\label{eq:proof-4}
    \begin{split}
        0 \leq& \displaystyle  ( \bm{u} - \bm{u}_{t+1})^{\top} \bm{\xi}_{t} + \frac{1}{2\gamma} (\| \bm{u} -\bm{u}_{t}\|_2^2 -\| \bm{u} - \bm{u}_{t+1}\|_2^2 - \| \bm{u} _{t+1} - \bm{u}_{t}\|_2^2) \\ 
        & \displaystyle + \dfrac{\sigma}{2} (\| \bm{u} - \hat{\bm{u}}\|_2 ^2 -\| \bm{u}- \bm{u}_{t+1}\|_2 ^2- \| \bm{u} _{t+1} - \hat{\bm{u}}\|_2 ^2) \quad\text{and} \\
        0 \leq& \displaystyle   (\bm{v}- \bm{v}_{t+1})^{\top} \left(2\bm{\zeta}_{t}    - \bm{\zeta}_{t-1}  \right)+ 
        \frac{1}{2\delta }\|\bm{v} - \bm{v}_{t}\|_2^2 - \frac{1}{2\delta } \| \bm{v}- \bm{v}_{t+1}\|_2^2 -  \frac{1}{2\delta } \| \bm{v}_{t+1}, \bm{v}_t\|_2^2.
    \end{split}
\end{equation}
Combining \eqref{equality: thm-iter-comp-algo2-1}, \eqref{eq:proof-5} and \eqref{eq:proof-4} yields 
\begin{equation}\label{eq:proof-7}
    \begin{split}
        \hat{\mathcal{F}}(\bm{u} _{t+1}, \bm{v}) - \hat{\mathcal{F}}( \bm{u}, \bm{v}_{t+1})
        \leq &\,\frac{1}{2\gamma}\|\bm{u} -\bm{u}_{t}\|_2^2  -\frac{1}{2\gamma}\|\bm{u}- \bm{u}_{t+1}\|_2^2 - \frac{1}{2\gamma}\|\bm{u}_{t+1} - \bm{u}_{t}\|_2^2  \\
        & + \frac{1}{2\delta} \|\bm{v} -\bm{v}_{t}\|_2^2 -\frac{1}{2\delta} \| \bm{v}- \bm{v}_{t+1}\|_2^2 - \frac{1}{2\delta}   \|\bm{v}_{t+1}- \bm{v}_t\|_2^2\\
        &  + \dfrac{\sigma}{2} \|\bm{u} - \hat{\bm{u}}\|_2^2 -\frac{\sigma}{2} \|\bm{u}- \bm{u}_{t+1}\|_2 ^2  - \frac{\sigma}{2}\|\bm{u}_{t+1} -\hat{\bm{u}}\|_2^2  \\
        & + (\bm{v}- \bm{v}_{t+1})^{\top} \left(2\bm{\zeta}_{t} - \bm{\zeta}_{t+1} - \bm{\zeta}_{t-1}  \right) + (\bm{u}_{t} - \bm{u}_{t+1})^{\top}(\bm{\xi}_{t}- \bm{\xi}'_{t}) .
    \end{split}
\end{equation}
Using the bounds on the subgradients and the fact
\[
(\bm{v}_{t+1} -\bm{v}_t)^{\top}(\bm{\zeta}_{t-1}- \bm{\zeta}_t) \leq \frac{\delta}{2} \|\bm{\zeta}_{t-1}-\bm{\zeta}_{t} \|_2^2 + \frac{1}{2\delta} \| \bm{v}_{t+1}- \bm{v}_{t}\|_2^2 \leq 2\delta C_2^2 + \frac{1}{2\delta} \| \bm{v}_{t+1}- \bm{v}_{t}\|_2^2,
\]
we obtain the two inequalities
\begin{equation}\label{eq:proof-6}
    \begin{split}
        (\bm{u} _{t} - \bm{u}_{t+1})^{\top}(\bm{\xi}_{t}- \bm{\xi}'_{t}) &\leq \displaystyle \frac{1}{2\gamma} \|\bm{u}_t -\bm{u}_{t+1}\|_2^2 + 2\gamma C_1^2 ~~~{\rm and}~~~ \\
        (\bm{v}- \bm{v}_{t+1})^{\top} \left(2\bm{\zeta}_{t} - \bm{\zeta}_{t+1} - \bm{\zeta}_{t-1} \right) &\leq  \displaystyle (\bm{v}- \bm{v}_{t+1})^{\top}(\bm{\zeta}_t- \bm{\zeta}_{t+1}) - 
        (\bm{v}- \bm{v}_{t})^{\top}(\bm{\zeta}_{t-1}- \bm{\zeta}_t)  \\
        & \quad +  2\delta C_2^2 + \frac{1}{2\delta} \| \bm{v}_{t+1} - \bm{v}_{t} \|_2^2,
    \end{split}
\end{equation}
Substituting \eqref{eq:proof-6} into \eqref{eq:proof-7}, we get
\begin{align*}
    &\, \hat{\mathcal{F}}( \bm{u} _{t+1}, \bm{v}) - \hat{\mathcal{F}}( \bm{u}, \bm{v}_{t+1}) + \dfrac{\sigma}{2} ( \| \bm{u} _{t+1} -\hat{\bm{u}}\|_2 ^2 + \| \bm{u}_{t+1} -\bm{u} \| _2 ^2) \\
    \le & \, \dfrac{\sigma}{2} \| \bm{u} - \hat{\bm{u}}\|_2^2 + \frac{1}{2\gamma} \|\bm{u} -\bm{u}_{t}\|_2^2 
    + \frac{1}{2\delta}  \| \bm{v}- \bm{v}_{t}\|_2^2  - \frac{1}{2\gamma}\|\bm{u} - \bm{u}_{t+1}\|_2^2-\frac{1}{2\delta}  \| \bm{v}- \bm{v}_{t+1}\|_2^2\\
    & + (\bm{v}- \bm{v}_{t+1})^{\top}(\bm{\zeta}_t- \bm{\zeta}_{t+1}) - 
    (\bm{v}- \bm{v}_{t})^{\top}(\bm{\zeta}_{t-1}- \bm{\zeta}_t) + 2\gamma C_1^2 +  2\delta C_2^2.
\end{align*}
Summing over $t=0, \cdots, T-1$ and using the convexity of $\hat{\mathcal{F}}(\cdot, \bm{v})$, $- \hat{\mathcal{F}}( \bm{u}, \cdot)$ and $\|\cdot\|_2 ^2$, we have that for any $(\bm{u}, \bm{v}) \in \mathcal{U} \times \mathcal{V}$,
\begin{align*}
    \hat{\mathcal{F}}(\bar{\bm{u}}_T, \bm{v}) - \hat{\mathcal{F}}(\bm{u}, \bar{\bm{v}}_T) 
    &\le  \frac{\|\bm{u} - \bm{u}_{0} \|_2^2}{2\gamma T} + \frac{\|\bm{v} - \bm{v}_0 \|_2^2}{2\delta T} + 2  C_1^2 \gamma+ 2  C_2^2 \delta \\
    & \quad + \dfrac{\sigma}{2} (\|\bm{u} - \hat{\bm{u}}\|_2^2- \|\bar{\bm{u}}_T- \hat{\bm{u}}\|_2^2 - \|\bm{u}- \bar{\bm{u}}_T\|_2^2) \\
    & \le   \mathcal{O}(T^{-1/2})  + \dfrac{\sigma}{2} (\|\bm{u} - \hat{\bm{u}}\|_2^2- \|\bar{\bm{u}}_T - \hat{\bm{u}}\|_2^2 - \|\bm{u}- \bar{\bm{u}}_T \|_2^2),
\end{align*}
where the first inequality follows from the definitions of $\bar{\bm{u}}_T$ and $\bar{\bm{v}}_T$, the bounds on subgradients, and the inequality
\begin{equation}
    \label{eq:proof-8}
    (\bm{v} -\bm{v}_{T})^{\top}(\bm{\zeta}_{T-1}- \bm{\zeta}_{T}) \leq \frac{\delta}{2} \|\bm{\zeta}_{T-1}-\bm{\zeta}_{T} \|_2^2 + \frac{1}{2\delta} \|\bm{v}- \bm{v}_{T} \|_2^2\leq  2C_2^2 \delta + \frac{1}{2\delta} \|\bm{v}- \bm{v}_{T} \|_2^2,
\end{equation}
and the second inequality follows from $\gamma = \delta = T^{-1/2}$. This completes the proof.
\end{proof}

\begin{proposition}\label{thm: iteration-complexity-Algo2-smooth}
    Instate the conditions of Proposition~\ref{thm: iteration-complexity-Algo2}.
    Suppose furthermore the existence of $\mathrm{Lip}_{\bm{u}\bm{u}}, \mathrm{Lip}_{\bm{v}\bm{v}}, \mathrm{Lip}_{\bm{v}\bm{u}}>0$ such that for any $\bm{u},\bm{u}'\in \mathcal{U}$, $\bm{v},\bm{v}'\in \mathcal{V}$, 
    \begin{equation}\label{Assump:u} 
    \begin{split}
        \| \nabla_{\bm{u}}\hat{\mathcal{F}}(\bm{u}, \bm{v}) - \nabla_{\bm{u}}\hat{\mathcal{F}}(\bm{u}', \bm{v})\|_2 &\leq \mathrm{Lip}_{\bm{u}\bm{u}} \|\bm{u} - \bm{u}'\|_2 \quad\text{and} \\
        \|\nabla_{\bm{v}}\hat{\mathcal{F}}(\bm{u}, \bm{v}) - \nabla_{\bm{v}}\hat{\mathcal{F}}(\bm{u}', \bmt{v}) \|_2 &\leq \mathrm{Lip}_{\bm{v}\bm{v}} \|\bm{v} - \bm{v}'\|_2+ \mathrm{Lip}_{\bm{v}\bm{u}} \|\bm{u} - \bm{u}'\|_2.
    \end{split}
    \end{equation}
    Then, the output $(\bar{\bm{u}}_T, \bar{\bm{v}}_T)$ of Algorithm~\ref{alg:mpsad} with $\gamma\le    \frac{1}{\mathrm{Lip}_{\bm{u}\bm{u}}+\mathrm{Lip}_{\bm{v}\bm{u}}}$ and $\delta \le \frac{1}{ 2\mathrm{Lip}_{\bm{v}\bm{v}}+\mathrm{Lip}_{\bm{v}\bm{u}}}$ is a strong $\mathcal{O}(T^{-1})$-approximate saddle point of problem~\eqref{opt:inner-saddle-general}. 
\end{proposition}

\begin{proof}[Proof of \Cref{thm: iteration-complexity-Algo2-smooth}.]
Using the Lipschitz property of the gradient of $\hat{\mathcal{F}}$, we can improve the second inequality in \eqref{eq:proof-5} to
\begin{equation*}
    \hat{\mathcal{F}}(\bm{u}_{t+1}, \bm{v}_{t+1}) - \hat{\mathcal{F}}(\bm{u}_t, \bm{v}_{t+1}) \le  (\bm{u}_{t+1} - \bm{u}_{t} )^{\top} \nabla_{\bm{u}}\hat{\mathcal{F}}(\bm{u}_t, \bm{v}_{t+1}) + \frac{\mathrm{Lip}_{\bm{u}\bm{u}}}{2} \|\bm{u}_t - \bm{u}_{t+1}\|_2^2,
\end{equation*}
the two inequalities in \eqref{eq:proof-6} to
\begin{align*}
    (\bm{u} _{t} - \bm{u}_{t+1})^{\top}(\bm{\xi}_{t}- \bm{\xi}'_{t})& =0 \quad\text{and}\\
    (\bm{v}- \bm{v}_{t+1})^{\top} \left( 2\bm{\zeta}_{t} - \bm{\zeta}_{t+1} - \bm{\zeta}_{t-1} \right)
    &\leq\displaystyle (\bm{v}- \bm{v}_{t+1})^{\top}(\bm{\zeta}_t- \bm{\zeta}_{t+1}) - (\bm{v}- \bm{v}_{t})^{\top}(\bm{\zeta}_{t-1}- \bm{\zeta}_t)\\
    & +\frac{\mathrm{Lip}_{\bm{v}\bm{v}}+ \mathrm{Lip}_{\bm{v}\bm{u}}}{2}\|\bm{v}_t - \bm{v}_{t+1}\|_2^2 + \frac{\mathrm{Lip}_{\bm{v}\bm{v}}}{2} \|\bm{v}_t - \bm{v}_{t-1}\|_2^2 + \frac{\mathrm{Lip}_{\bm{v}\bm{u}}}{2}\|\bm{u}_t - \bm{u}_{t-1}\|_2^2,
\end{align*}
and inequality~\eqref{eq:proof-8} to
\begin{equation*}
    (\bm{v}- \bm{v}_{T})^{\top}(\bm{\zeta}_{T-1} - \bm{\zeta}_{T}) \leq  \displaystyle \frac{\mathrm{Lip}_{\bm{v}\bm{v}}+ \mathrm{Lip}_{\bm{v}\bm{u}}}{2}\|\bm{v} - \bm{v}_{T}\|_2^2 + \frac{\mathrm{Lip}_{\bm{v}\bm{v}}}{2} \|\bm{v}_{T} - \bm{v}_{T-1}\|_2^2 + \frac{\mathrm{Lip}_{\bm{v}\bm{u}}}{2}\|\bm{u}_{T} - \bm{u}_{T-1}\|_2^2 .
\end{equation*}
Following the proof of Proposition~\ref{thm: iteration-complexity-Algo2} but using the above improved inequalities, we get
\begin{align*}
    &\,\hat{\mathcal{F}}(\bar{\bm{u}}_T, \bm{v}) - \hat{\mathcal{F}}(\bm{u}, \bar{\bm{v}}_T) + \dfrac{ \sigma}{2} \left( \| \bar{\bm{u}}_T -\hat{\bm{u}}\|_2^2 +  \|\bm{u}- \bar{\bm{u}}_T\|_2 ^2  - \|\bm{u} - \hat{\bm{u}}\|_2^2 \right) \\
    \le&\, \frac{1}{2T\gamma} \|\bm{u} -\bm{u}_{0}\|_2^2  +  \frac{1}{2\delta T} \|\bm{v} -  \bm{v}_{0} \|_2^2 -\frac{1}{2T\gamma} \|\bm{u}- \bm{u}_{T}\|_2^2 -\frac{1}{T} \left(\frac{1}{2\gamma} - \frac{\mathrm{Lip}_{\bm{u}\bm{u}} +\mathrm{Lip}_{\bm{v}\bm{u}}}{2} \right) \sum_{t=1}^{T} \|\bm{u}_{t} - \bm{u}_{t-1}\|_2^2   \\
    &   -\frac{1}{T} \left(\frac{1}{2 \delta}- \frac{2\mathrm{Lip}_{\bm{v}\bm{v}}+ \mathrm{Lip}_{\bm{v}\bm{u}}}{2} \right) \sum_{t=1}^{T}\| \bm{v}_{t} - \bm{v}_{t-1}\|_2^2 - \frac{1}{T} \left(\frac{1}{2  \delta}- \frac{\mathrm{Lip}_{\bm{v}\bm{v}}+ \mathrm{Lip}_{\bm{v}\bm{u}}}{2} \right) \|\bm{v}_{T-1}-\bm{v}_T\|_2^2 \\
    \le&\,  \mathcal{O}(T^{-1}), 
\end{align*}
where the second inequality follows from $\gamma\le    \frac{1}{\mathrm{Lip}_{\bm{u}\bm{u}}+\mathrm{Lip}_{\bm{v}\bm{u}}}$ and $\delta \le \frac{1}{ 2\mathrm{Lip}_{\bm{v}\bm{v}}+\mathrm{Lip}_{\bm{v}\bm{u}}}$.
\end{proof}

We are now ready to prove \Cref{thm:inner-rate}.

\begin{proof}[Proof of \Cref{thm:inner-rate}.]
The \ref{opt:inner-saddle} problem is a special case of problem~\eqref{opt:inner-saddle-general}, with $\sigma = 1/\alpha$, $\mathcal{U} = \mathcal{X}$, $\mathcal{V} = \mathcal{Z}$ and
\begin{equation}
    \label{eq:proof-9}
    \hat{\mathcal{F}} (\bm{x}, \bm{z}) = f_0(\bm{x}) + \sum_{m\in [M] }\lambda_{m}^{k+1} g_m(\bm{x}, \bm{z}_m) .
\end{equation}
When applied to the \ref{opt:inner-saddle} problem, Algorithm~\ref{alg:inner} reduces to Algorithm~\ref{alg:mpsad}.
By Assumption~\ref{Assump: Robust-Prob-1}, we see that $\hat{\mathcal{F}}(\cdot, \bm{z})$ is convex on $\mathcal{X}$ for any $\bm{z}\in \mathcal{Z}$ and $\hat{\mathcal{F}}(\bm{x}, \cdot)$ is concave on $\mathcal{Z}$ for any $\bm{x}\in\mathcal{X}$. Also, by Assumption~\ref{Assump: Robust-Prob-2}, we have that the subdifferentials $\partial_{\bm{x}} \hat{\mathcal{F}}(\bm{x}, \bm{z})$ and $\partial_{\bm{z}} (-\hat{\mathcal{F}})(\bm{x}, \bm{z})$ are both uniformly bounded on $\mathcal{X}\times \mathcal{Z}$.
The desired result thus follows from \Cref{thm: iteration-complexity-Algo2}.
\end{proof}

\Cref{thm:inner-rate-smooth} can be proved in a similar vein.

\begin{proof}[Proof of \Cref{thm:inner-rate-smooth}.]
The proof is similar to that of \Cref{thm:inner-rate}, except that we need to verify the Lipschitz property \eqref{Assump:u} and determine the step-sizes.
Recall that for the \ref{opt:inner-saddle} problem, $\hat{\mathcal{F}}(\bm{x}, \bm{z})$ is given by \eqref{eq:proof-9}. Therefore, by Assumption~\ref{ass:smoothness},
\begin{align*}
    &\, \| \nabla_{\bm{x}} \hat{\mathcal{F}}(\bm{x}, \bm{z}) - \nabla_{\bm{x}} \hat{\mathcal{F}}(\bm{x}', \bm{z}) \|_2 \\
    =&\, \left\| \nabla_{\bm{x}} f_0(\bm{x}) - \nabla_{\bm{x}} f_0(\bm{x}') \right\|_2 + \sum_{m\in [M]} \lambda_m^{k+1} \| \nabla_{\bm{x}} g_m (\bm{x}, \bm{z}_m) - \nabla_{\bm{x}} g_m (\bm{x}', \bm{z}_m) \|_2\\
    \le &\, \left(D_0'  + \sum_{m\in [M]} \lambda_m^{k+1} D'_m \right) \| \bm{x} - \bm{x}' \|_2,
\end{align*}
which implies $\mathrm{Lip}_{\bm{u}\bm{u}} = (D_0'  + \sum_{m\in [M]} \lambda_m^{k+1} D'_m )$.

Next, by Assumption~\ref{ass:smoothness},
\begin{align*}
    &\, \| \nabla_{\bm{z}} \hat{\mathcal{F}}(\bm{x}, \bm{z}) - \nabla_{\bm{z}} \hat{\mathcal{F}}(\bm{x}', \bm{z}') \|_2 \\
    =&\, \sqrt{ \sum_{m\in[M]}  {\lambda_m^{k+1}}^2 \| \nabla_{\bm{z}_m} g_m (\bm{x}, \bm{z}_m) - \nabla_{\bm{z}_m} g_m (\bm{x}', \bm{z}_m') \|_2^2} \\
    \le &\, \sqrt{ \sum_{m\in[M]}  {\lambda_m^{k+1}}^2 \left( E'_{m,1}\| \bm{x} - \bm{x}' \|_2 + E'_{m,2} \| \bm{z}_m - \bm{z}_m' \|_2 \right)^2 } \\
    \le &\, \sqrt{2 \sum_{m\in[M]} {\lambda_m^{k+1}}^2 {E'_{m,1}}^2 } \| \bm{x} - \bm{x}' \|_2 + \sqrt{ 2\max_{m\in[M]} {\lambda_m^{k+1}}^2 {E'_{m,2}}^2 } \|\bm{z} - \bm{z}'\|_2,
\end{align*}
which implies $\mathrm{Lip}_{\bm{v}\bm{u}} = \sqrt{2 \sum_{m\in[M]} {\lambda_m^{k+1}}^2 {E'_{m,1}}^2 }$ and $\mathrm{Lip}_{\bm{v}\bm{v}} = \sqrt{ 2\max_{m\in[M]} {\lambda_m^{k+1}}^2 {E'_{m,2}}^2 }$.
The desired result thus follows from \Cref{thm: iteration-complexity-Algo2-smooth}.
\end{proof}

\section{Oracle Complexity}

\begin{proof}[Proof of \Cref{thm:complexity_1}.]
Note that the outer algorithm does not directly rely on the projection or subgradient oracles, but only through the invocation of the inner algorithm. Also, each iteration of the inner algorithm requires at most a constant number of calls to the oracles, independent of $\varepsilon$. Therefore, to prove the oracle complexity of \prom, it suffices to count the aggregated number of inner iterations.

Consider Algorithms~\ref{alg:outer} and~\ref{alg:inner} with $K = \mathcal{O}(\varepsilon^{-1})$ and other algorithmic parameters chosen as in Theorems~\ref{thm: iteration-complexity-Algo1} and~\ref{thm:inner-rate}. 
By \Cref{thm: iteration-complexity-Algo1}, \prom produces an $\varepsilon$-approximate optimal solution to the \ref{opt:RO1} problem in $\mathcal{O}(\varepsilon^{-1})$ outer iterations. So, we need to invoke the inner algorithm $\mathcal{O}(\varepsilon^{-1})$ times. Each invocation needs to compute a strong $ \mathcal{\varepsilon}$-approximate saddle point, which by \Cref{thm:inner-rate}, requires $\mathcal{O}(\varepsilon^{-2})$ inner iterations. We thus conclude that the oracle complexity is $\mathcal{O}(\varepsilon^{-1})\mathcal{O}(\varepsilon^{-2}) = \mathcal{O}(\varepsilon^{-3})$.
\end{proof}

\begin{proof}[Proof of \Cref{thm:complexity_2}.]
The proof is the same as that of \Cref{thm:complexity_1} (except that we use \Cref{thm:inner-rate-smooth} instead of \Cref{thm:inner-rate}) and hence omitted.
\end{proof}

\section{Proofs for Extended \prom}\label{sec:extended_prom}

\begin{proof}[Proof of Proposition~\ref{prop:extended_RO}]
    Consider the \ref{opt:RO1} problem with projection-unfriendly uncertainty sets~\eqref{def:generalized-Zm}. For any fixed $m\in[M]$ and $\bm{x}\in\mathcal{X}$, the embedded maximization problem in the robust constraint reads
    \begin{equation}
        \label{opt:embed_extend}
            \begin{array}{rll}
            \max & g_m (\bm{x}, \bm{z}_m)  \\
            {\rm s.t.} & h_{m,i} ( \bm{z}_m) \leq 0  \quad\forall i\in [I_m] \\
            &  \bm{z}_m \in \widetilde{\mathcal{Z}}_m.
            \end{array}
    \end{equation}
    By Assumption~\ref{Assump: Robust-Prob-1-gen}\ref{Assump: Robust-Prob-1-gen-i}-\ref{Assump: Robust-Prob-1-gen-iii}, we have that $\widetilde{\mathcal{Z}}_m$ is a non-empty, compact and convex set, that $g_m(\bm{x}, \cdot)$ is concave and $h_{m,i}$ is convex for all $i\in [I_m]$, and that a Slater point for problem~\eqref{opt:embed_extend} exists. Therefore, strong duality holds and problem~\eqref{opt:embed_extend} is equivalent to its Lagrangian dual
    \begin{equation}
        \label{opt:embed_extend_dual}
            \min_{\bm{\mu}_m \in\mathbb{R}_+^{I_m}} \max_{ \bm{z}_m\in \widetilde{\mathcal{Z}}_m }\; g_m (\bm{x}, \bm{z}_m) - \bm{\mu}_m^\top \bm{h}_m (\bm{z}_m)
    \end{equation}
    Replacing the embedded problem by problem~\eqref{opt:embed_extend_dual}, we see that the \ref{opt:RO1} problem is equivalent to 
    \begin{equation*}
        \begin{array}{rll}
                \min & f_0(\bm{x}) \\
                {\rm s.t.} & \displaystyle  \max_{\bm{z}_m \in \widetilde{\mathcal{Z}}_m} g_m (\bm{x}, \bm{z}_m) - \bm{\mu}_m^\top \bm{h}_m (\bm{z}_m) \leq 0  \quad\forall m\in [M] \\
                & \bm{x} \in \mathcal{X},\; \bm{\mu}_m  \in \mathbb{R}_+^{I_m}\quad\forall m\in [M].
            \end{array}
    \end{equation*}
    It remains to prove that any optimal solution $(\bm{x}^\star, \bm{\mu}^\star) \in\mathcal{X}\times \mathbb{R}_+^{I_1+\cdots+I_M}$ satisfies that
    \begin{equation*}
        \mu_{m,i}^\star \le \frac{G_m}{\max_{i\in [I_{m}]} \{h_{m, i}(\bar{\bm{z}}_m)\}}\quad\forall i\in [I_m], \, m\in [M].
    \end{equation*}
    To do so, let $(\bm{x}^\star, \bm{\mu}^\star) \in\mathcal{X}\times \mathbb{R}_+^{I_1+\cdots+I_M}$ be any optimal solution. Then, for any $m\in [M]$,
    \begin{align*}
        \max_{\bm{z}_m \in \widetilde{\mathcal{Z}}_m} g_m (\bm{x}^\star, \bm{z}_m) - {\bm{\mu}_m^\star}^\top \bm{h}_m (\bm{z}_m) \leq 0,
    \end{align*}
    which implies that
    \begin{equation*}
        g_m (\bm{x}^\star, \bar{\bm{z}}_m) \le  {\bm{\mu}_m^\star}^\top \bm{h}_m (\bar{\bm{z}}_m) .
    \end{equation*}
    Noting that $\mu_{m,i}^\star \, h_{m,i}(\bar{\bm{z}}_m) \le 0$ for all $i\in [I_m]$, we have 
    \begin{equation*}
        \mu_{m,i}^\star \, h_{m,i}(\bar{\bm{z}}_m) \ge {\bm{\mu}_m^\star}^\top \bm{h}_m (\bar{\bm{z}}_m) \ge  g_m (\bm{x}^\star, \bar{\bm{z}}_m) \ge G_m.
    \end{equation*}
    By Assumption~\ref{Assump: Robust-Prob-1-gen}\ref{Assump: Robust-Prob-1-gen-iii}, $h_{m,i}(\bar{\bm{z}}_m) < 0$ for all $i\in [I_m]$. Therefore, for any $i\in [I_m]$,
    \begin{equation*}
        \mu^\star_{m,i} \le \frac{G_m}{h_{m,i}(\bar{\bm{z}}_m)} \le \frac{G_m}{\max_{i\in[I_m]}\{h_{m,i}(\bar{\bm{z}}_m)\}}.
    \end{equation*}
    This completes the proof.
\end{proof}

\begin{proof}[Proof of Theorem~\ref{thm:extended_prom_complexity}]
    The extended \prom is the algorithm obtained by applying \prom (from Section~\ref{sec:prom}) to the problem~\ref{opt:Gene-RO1}.
    It suffices is to verify that problem~\ref{opt:Gene-RO1} satisfies Assumptions~\ref{Assump: Robust-Prob-1}, \ref{Assump: Robust-Prob-2} and~\ref{ass:smoothness} in the sense that these assumptions hold when the data $(f_0,g_1,\dots, g_M, \mathcal{X}, \mathcal{Z}_1,\dots, \mathcal{Z}_M)$ is replaced by $ (\tilde{f}_0, \tilde{g}_1,\dots, \tilde{g}_M, \widetilde{\mathcal{X}}, \widetilde{\mathcal{Z}}_1,\dots, \widetilde{\mathcal{Z}}_M) $.

    \textbf{Assumption~\ref{Assump: Robust-Prob-1}\ref{Assump: Robust-Prob-1-i}:} The non-emptiness, compactness and convexity of the sets $\widetilde{\mathcal{Z}}_1,\dots,\widetilde{\mathcal{Z}}_M$ follow directly from Assumption~\ref{Assump: Robust-Prob-1-gen}\ref{Assump: Robust-Prob-1-gen-i}. Recall that $\widetilde{\mathcal{X}} = \mathcal{X} \times \mathcal{M}$, where $\mathcal{M} = [0, a_1]^{I_1}\times\cdots\times [0, a_M]^{I_M}$. Assumption~\ref{Assump: Robust-Prob-1-gen}\ref{Assump: Robust-Prob-1-gen-i} therefore implies also that $\widetilde{\mathcal{X}}$ is non-empty, compact and convex.
    
    \textbf{Assumption~\ref{Assump: Robust-Prob-1}\ref{Assump: Robust-Prob-1-ii}:} The objective function $\tilde{f}_0(\tilde{\bm{x}}) = f_0 (\bm{x})$ is obviously convex and finite-valued on $\widetilde{\mathcal{X}} = \mathcal{X}\times \mathcal{M}$, since $f_0$ is convex and finite-valued on $\mathcal{X}$ by Assumption~\ref{Assump: Robust-Prob-1-gen}\ref{Assump: Robust-Prob-1-gen-ii}. Also by Assumption~\ref{Assump: Robust-Prob-1-gen}\ref{Assump: Robust-Prob-1-gen-ii}, $g_m(\cdot, \bm{z}_m)$ is convex for any $\bm{z}_m\in \widetilde{\mathcal{Z}}_m$ and $g_m (\bm{x}, \cdot) $ is concave for any $\bm{x}\in \mathcal{X}$, and $h_{m,i}$ is convex and finite-valued on $\widetilde{\mathcal{Z}}_m$ for all $i\in [I_m]$ and $m\in[M]$. Together with the boundedness of $\mathcal{M}$, this implies that the function $\tilde{g}_m(\bmt{x}, \bm{z}_m) = g_m (\bm{x}, \bm{z}_m)- \bm{\mu}_m^\top \bm{h}_m (\bm{z}_m) $ is finite-valued on $\widetilde{\mathcal{X}}\times \widetilde{\mathcal{Z}}_m$ and satisfies that $\tilde{g}_m (\cdot, \bm{z}_m)$ is convex for any $\bm{z}\in \widetilde{\mathcal{Z}}_m$ and $\tilde{g}_m(\bmt{x}, \cdot)$ is concave for any $\bmt{x} = (\bm{x}, \bm{\mu}) \in \widetilde{\mathcal{X}}$.

    \textbf{Assumption~\ref{Assump: Robust-Prob-1}\ref{Assump: Robust-Prob-1-iii}:} By Assumption~\ref{Assump: Robust-Prob-1-gen}\ref{Assump: Robust-Prob-1-gen-iii}, there exists $\bar{\bm{x}}$ such that $\displaystyle \max_{\bm{z}_m \in \widetilde{\mathcal{Z}}_m} g_m (\bar{\bm{x}}, \bm{z}_m) <0$ for any $m\in[M]$. Therefore, $\max_{\bm{z}_m \in \widetilde{\mathcal{Z}}_m} \tilde{g}_m( (\bar{\bm{x}}, \bm{0}) , \bm{z}_m) = \max_{\bm{z}_m \in \widetilde{\mathcal{Z}}_m} g_m (\bar{\bm{x}}, \bm{z}_m) <0$.

    \textbf{Assumption~\ref{Assump: Robust-Prob-1}\ref{Assump: Robust-Prob-1-iv}:} It follows directly from Assumption~\ref{Assump: Robust-Prob-1-gen}\ref{Assump: Robust-Prob-1-gen-iv}.

    \textbf{Assumption~\ref{Assump: Robust-Prob-2}:} By Assumption~\ref{Assump: Robust-Prob-2-gen}, we have that $\tilde{f}_0 (\tilde{\bm{x}}) = f_0 (\bm{x})$ is subdifferentiable on $\widetilde{\mathcal{X}} = \mathcal{X}\times\mathcal{M}$. Recall that $\tilde{g}_m (\tilde{\bm{x}}, \bm{z}_m) = g_m (\bm{x}, \bm{z}_m)- \bm{\mu}_m^\top \bm{h}_m (\bm{z}_m) $ for any $m\in[M]$. By Assumption~\ref{Assump: Robust-Prob-2-gen}, we have that $\tilde{g}_m (\tilde{\bm{x}}, \bm{z}_m)  $ is subdifferentiable in $\bmt{x} = (\bm{x}, \bm{\mu})$ on $ \widetilde{\mathcal{X}} $ and $- \tilde{g}_m (\tilde{\bm{x}}, \bm{z}_m)  $ subdifferentiable in $\bm{z}_m$ on $ \widetilde{\mathcal{Z}}_m$. 
    
    We then bound the subdifferentials. First, any subgradient of $\tilde{f}_0 $ at $\bmt{x}\in\widetilde{\mathcal{X}}$ is of the form $\tilde{\bm{\xi}}_0 = (\bm{\xi}_0, \bm{0})$ for some $\bm{\xi}_0\in \partial f_0 (\bm{x})$. So, by Assumption~\ref{Assump: Robust-Prob-2-gen}, the subdifferential $\partial \tilde{f}_0 (\bmt{x})$ is uniformly bounded by $\tilde{D}_0 = D_0$. 
    Next, for any $m\in [M]$ and $\bm{z}_m\in  \widetilde{\mathcal{Z}}_m$, any subgradient of $\tilde{g}_m (\cdot, \bm{z}_m)$ is of the form $\tilde{\bm{\xi}}_m = ( \bm{\xi}_m, \bm{0}, -\bm{h}_m (\bm{z}_m), \bm{0} )$ for some $ \bm{\xi}_m\in \partial_{\bm{x}} g_m(\bm{x}, \bm{z}_m)$.
    Since $H_m$ is real-valued and convex on $\widetilde{\mathcal{Z}}_m$, there exists a constant $H_m > 0$ such that $\|\bm{h}_m( \bm{z}_m)\| \le H_m$ for all $\bm{z}_m \in \mathcal{Z}_m$.
    Noting that $\| \tilde{\bm{\xi}}_m \| \le \sqrt{ \| \bm{\xi}_m \|^2 + \| \bm{h}_m(\bm{z}_m) \|^2 }$, by Assumption~\ref{Assump: Robust-Prob-2-gen}, the subdifferential $ \partial_{\bmt{x}} \tilde{g}_m (\bmt{x}, \bm{z}_m) $ is uniformly bounded by $\tilde{D}_m = \sqrt{D_m^2 + H_m^2}$. 
    Finally, for any $m\in [M]$ and $\bmt{x} = (\bm{x}, \bm{\mu})\in \widetilde{\mathcal{X}} $, any subgradient of $(-\tilde{g}_m) (\bmt{x}, \cdot)$ is of the form $\tilde{\bm{\zeta}}_m = \bm{\zeta}_m - \sum_{i\in [I_m]} \mu_{m,i}\bm{\eta}_{m,i} $ for some $\bm{\eta}_{m,1} \in \partial h_{m,1} (\bm{z}_m),\dots, \bm{\eta}_{m, I_m} \in \partial h_{m,I_m} (\bm{z}_m)$ and $ \bm{\zeta}_m\in \partial_{\bm{z}_m} (-g_m)(\bm{x}, \bm{z}_m)$. 
    Noting that $\| \tilde{\bm{\zeta}}_m \| \le \| \bm{\zeta}_m \| + \sum_{i\in [I_m]} \mu_{m,i} \| \bm{\eta}_{m,i}\|$, by Assumption~\ref{Assump: Robust-Prob-2-gen}, the subdifferential $ \partial_{\bm{z}_m} (-\tilde{g}_m) (\bmt{x}, \bm{z}_m) $ is uniformly bounded by $\tilde{E}_m = E_m + a_m I_m F_m$.

    \textbf{Assumption~\ref{ass:smoothness}:} Since $\tilde{f}_0 (\tilde{\bm{x}}) = f_0 (\bm{x})$ for any $\bmt{x} \in \widetilde{\mathcal{X}}$, it follows from Assumption~\ref{ass:smoothness-gen} that $\nabla \tilde{f}_0$ is Lipschitz continuous on $\widetilde{\mathcal{X}}$. For any $m\in[M]$, $\bm{z}\in \widetilde{\mathcal{Z}}_m$ and $\bmt{x} = (\bm{x}, \bm{\mu})\in \widetilde{\mathcal{X}}$, $\nabla_{\bmt{x}} \tilde{g}_m (\bmt{x}, \bm{z}_m) = (\nabla_{\bm{x}} g_m (\bm{x}, \bm{z}_m) ,\bm{0}, - \bm{h}_m (\bm{z}_m), \bm{0})$ and $\nabla_{\bm{z}_m} \tilde{g}_m (\bmt{x}, \bm{z}_m) = \nabla_{\bm{z}_m} g_m (\bm{x}, \bm{z}_m) -\sum_{i \in [I_m]} \mu_{m,i} \nabla h_{m,i} (\bm{z}_m)$. By Assumption~\ref{ass:smoothness-gen}, $\nabla_{\bmt{x}} \tilde{g}_m (\cdot, \bm{z}_m)$ is Lipschitz continuous on $\widetilde{\mathcal{X}}$ for any fixed $\bm{z}_m\in \widetilde{\mathcal{Z}}_m$, and $\nabla_{\bm{z}_m} \tilde{g}_m $ is jointly Lipschitz continuous on $\widetilde{X}\times \widetilde{\mathcal{Z}}_m$.

    The proof is completed.
\end{proof}

\end{document}